\documentclass[hidelinks,onefignum,onetabnum]{siamart251216}

\usepackage{comment}
\usepackage{amsfonts}
\usepackage{amssymb}
\usepackage{booktabs}
\usepackage{multirow}
\usepackage{graphicx}
\usepackage{xcolor}
\usepackage{tikz}
\usetikzlibrary{arrows.meta,calc,decorations.pathreplacing,positioning}
\usepackage{epstopdf}
\usepackage{algorithm}
\usepackage{algorithmic}
\ifpdf
  \DeclareGraphicsExtensions{.eps,.pdf,.png,.jpg}
\else
  \DeclareGraphicsExtensions{.eps}
\fi

\usepackage{enumitem}
\setlist[enumerate]{leftmargin=.5in}
\setlist[itemize]{leftmargin=.5in}

\crefname{theorem}{Theorem}{Theorems}
\newsiamthm{prediction}{Prediction}
\crefname{prediction}{Prediction}{Predictions}
\newsiamremark{remark}{Remark}
\crefname{remark}{Remark}{Remarks}
\newsiamremark{hypothesis}{Hypothesis}
\crefname{hypothesis}{Hypothesis}{Hypotheses}
\newsiamthm{claim}{Claim}
\newsiamremark{fact}{Fact}
\crefname{fact}{Fact}{Facts}

\headers{HiMuon: Tiled NS for Scalable LLM Pre-Training}{Authors}

\title{Hierarchical Muon: Tiled Newton--Schulz Updates for Efficient Muon Optimization}
\author{Ziyuan Tang\thanks{Department of Computer Science and Engineering, University of Minnesota, Minneapolis (\email{tang0389@umn.edu}, \email{saad@umn.edu}). The research of Tang and Saad is supported by the NSF award DMS-2513117.}
\and Tianshi Xu\thanks{Department of Mathematics, Emory University, Atlanta, GA 30322 (\email{tianshi.xu@emory.edu},
\email{yxi26@emory.edu}). The research of Xi and Xu is supported by NSF award  DMS-2513118.}
\and Yousef Saad\footnotemark[1]
\and Yuanzhe Xi\footnotemark[2]}

\headers{HiMuon}{Z. Tang, T. Xu, Y. Saad, and Y. Xi} 

\newcommand{\R}{\mathbb{R}}

\newcommand{\mA}{\mathbf{A}}

\newcommand{\mG}{\mathbf{G}}

\newcommand{\mI}{\mathbf{I}}
\newcommand{\mM}{\mathbf{M}}
\newcommand{\mU}{\mathbf{U}}
\newcommand{\mW}{\mathbf{W}}
\newcommand{\mX}{\mathbf{X}}
\newcommand{\mP}{\mathbf{P}}

\newcommand{\vtheta}{\boldsymbol{\theta}}

\DeclareMathOperator{\tr}{tr}

\DeclareMathOperator{\polar}{polar}

\newcommand{\Loss}{\mathcal{L}}

\newcommand{\normF}[1]{\left\lVert #1 \right\rVert_F}
\newcommand{\normTwo}[1]{\left\lVert #1 \right\rVert_2}

\newcommand{\lr}{\eta}
\newcommand{\tileSize}{T}

\newcommand{\nsA}{a}
\newcommand{\nsB}{b}
\newcommand{\nsC}{c}

\ifpdf
\hypersetup{
  pdftitle={Hierarchical Muon: A Tiled Newton--Schulz Iteration for Matrix-Function Preconditioning},
  pdfauthor={}
}
\fi

\begin{document}

\maketitle
\begin{abstract}
Muon-type optimizers construct update directions for dense neural-network
weights by applying a finite Newton--Schulz map to momentum-gradient matrices.
For an $H\times W$ matrix, with $r=\min\{H,W\}$ and $s=\max\{H,W\}$, $K$
steps of the full-matrix Newton--Schulz update require
$\mathcal{O}(r^2sK)$ work and couple all rows and columns through repeated
Gram matrix products.  We introduce Hierarchical Muon (HiMuon), a tiled
Newton--Schulz scheme for Muon-type optimization.  HiMuon partitions each
momentum-gradient matrix into $\tileSize\times\tileSize$ tiles, applies the
same finite Newton--Schulz map independently to each tile, and reassembles the
results.  For finite $\tileSize$ below the matrix dimensions, HiMuon defines a
local matrix-function map
rather than a convergent approximation to the full-matrix update: spectral
interactions are preserved within tiles and discarded across tile boundaries.
For fixed finite $\tileSize$, the leading Newton--Schulz work decreases to
$\mathcal{O}(HW\tileSize K)$, and the computation decomposes into independent
small dense matrix operations.  This structure enables
tile-size-dependent GPU kernels, cross-layer batching, memory-bounded
chunking, and runtime tile-size schedules.  Experiments on transformer
training and controlled matrix-function diagnostics show that HiMuon improves
optimizer-step efficiency while keeping training behavior close to
full-matrix Muon in the tested regimes.  Our implementation is publicly
available at \url{https://github.com/tang0389/himuon}.
\end{abstract}


\section{Introduction}
\label{sec:intro}

In transformer models, a large share of the parameter count and arithmetic
cost is concentrated in dense projection matrices in the attention and
feed-forward sublayers.  The gradients of these parameters inherit the same
matrix form.  Muon-type optimizers exploit this structure by applying a finite
Newton--Schulz (NS) iteration to the momentum-gradient matrix associated with
each matrix parameter~\cite{KellerJordan2024muon}.  The output is used as the
update direction.  This finite-step matrix function acts as a polar-factor
surrogate: it drives the singular values of the update direction toward a
common scale.  Such matrix-function updates have been shown empirically to
improve pretraining efficiency across several model scales~\cite{P1066,P1136,P1137}.

The full-matrix Newton--Schulz map is a major computational component of the
Muon optimizer step.  For an $H\times W$ momentum-gradient matrix, let
$r=\min\{H,W\}$ and $s=\max\{H,W\}$.  With the iteration oriented so that the
Gram matrix has dimension $r$, $K$ Newton--Schulz steps require
\[
    \mathcal{O}(r^2sK)
    =
    \mathcal{O}\!\left(HW\min\{H,W\}K\right)
\]
operations.  When $H\le W$, this becomes $\mathcal{O}(H^2WK)$.  The repeated
Gram matrix products also create global spectral coupling: an entry of the update can
depend on the entire momentum-gradient matrix.  Thus the computational cost of
the update and its spectral behavior are linked.  Retaining full-matrix
spectral coupling gives the standard Muon update, but doing so at every
optimizer step requires large dense matrix products. In distributed parameter-sharded training, the full-matrix Newton--Schulz iteration may also require reconstructing matrix shards through additional
collective communication inside the optimizer step.

In this paper, we study a tiled formulation in which the extent of this spectral
coupling is controlled explicitly by the tile size. Tiling simultaneously reduces the cost of dense matrix products, localizes the Newton--Schulz computation within shard-resident submatrices, and exposes more parallel work to the hardware. Full-matrix Muon
represents one extreme, where the Newton--Schulz map is applied to the entire
matrix.  Coordinatewise methods represent the opposite extreme, where no
matrix-level spectral interaction is retained.  Between these regimes is a
family of local matrix-function updates: the same Newton--Schulz map is
applied, but only to submatrices of a prescribed size.  The goal is to retain
useful local spectral behavior while reducing the cost of the optimizer step
and exposing more parallel work to the hardware.

We introduce Hierarchical Muon (HiMuon), a tiled Newton--Schulz method based on
this principle.  Let $\Phi_K$ denote the $K$-step Newton--Schulz map used in
Muon.  For a tile size $\tileSize$, HiMuon partitions a momentum-gradient
matrix into an $R\times C$ grid of $\tileSize\times\tileSize$ tiles, applies
$\Phi_K$ independently to each tile, and reassembles the tile outputs.  For
fixed finite $\tileSize$ below the matrix dimensions, this construction does
not recover the full-matrix Muon map.  It defines a distinct local
matrix-function update: spectral interactions are preserved within each tile
and removed across tile boundaries.  The tile size therefore acts as both a
numerical parameter, controlling the amount of retained spectral coupling, and
a computational parameter, controlling the cost of the Newton--Schulz
computation.  For fixed $\tileSize\ll\min\{H,W\}$, the leading
Newton--Schulz work scales as
\[
    \mathcal{O}(HW\tileSize K).
\]

The tiled formulation changes more than the arithmetic count.  Once the
partition is fixed, all tilewise Newton--Schulz evaluations are independent,
so the optimizer step becomes a batched small-matrix linear algebra workload
rather than a sequence of large full-matrix operations.  HiMuon exploits this
structure through tile-size-dependent kernels, cross-layer batching,
memory-bounded chunking, and runtime tile-size schedules.
The same representation also supports distributed layouts in which each rank
performs its tiled Newton--Schulz work locally.  In this sense, the numerical
locality of the update and the GPU execution strategy are designed together.

The main contributions of this work are summarized as follows.

\begin{itemize}
    \item \textbf{A tiled Newton--Schulz formulation of Muon.}
    We define HiMuon as a finite-tile matrix-function update for Muon-type
    optimization.  The formulation makes the degree of spectral coupling
    explicit through the tile size and distinguishes the full-matrix Muon
    update from local tiled alternatives.

    \item \textbf{An end-to-end GPU implementation.}
    We develop an implementation that maps the tiled Newton--Schulz workload to
    batched small dense matrix operations, using tile-size-dependent kernels,
    cross-layer batching, memory-bounded chunking, and support for distributed
    parameter layouts.

    \item \textbf{Empirical characterization across quality, speed, and precision.}
    We evaluate HiMuon at the matrix-function, layer, optimizer-step, and
    end-to-end training levels.  The experiments study real transformer weight
    shapes, kernel and tile-shape choices, cross-layer batching,
    finite-precision Newton--Schulz computation, and the resulting
    training-quality trade-off.  We also use an associative-memory benchmark as
    a structural diagnostic for tile-size selection.
\end{itemize}

The remainder of the paper is organized as follows.
\Cref{sec:prelim} introduces notation, the polar decomposition, the finite
Newton--Schulz iteration, and the baseline Muon optimizer.
\Cref{sec:related} reviews efficient variants of Muon-type optimization and
alternative reductions of the Newton--Schulz computation.
\Cref{sec:method} presents the HiMuon algorithm and its GPU implementation.
\Cref{sec:capacity} introduces a capacity-based diagnostic for the tiled map.
\Cref{sec:experiments} reports numerical experiments on training quality,
computational efficiency, cross-layer batching, and finite-precision behavior.
\Cref{sec:conclusion} concludes the paper.
\section{Preliminaries}
\label{sec:prelim}

This section fixes the notation and the full-matrix Muon update used as the baseline throughout the paper.  We first introduce matrix norms and the tiling
operator, then distinguish the exact polar factor from the finite
Newton--Schulz map used by Muon, and finally state the optimizer update that
HiMuon modifies.

\subsection{Notation}
\label{sec:notation}

Matrices are denoted by bold uppercase letters, e.g., $\mA$, and $\mI$ denotes
an identity matrix of the appropriate dimension.  The model parameters are
collected in $\vtheta$, and the training objective is $\Loss(\vtheta)$.
Subscript $t$ denotes the optimizer step.  For
$\mA\in\R^{H\times W}$, we write $\normF{\mA}$ and $\normTwo{\mA}$ for the
Frobenius and spectral norms.

For a finite tile size $\tileSize$, let $\mathcal{B}_{\tileSize}$ denote the
tiling operator.  Given $\mA\in\R^{H\times W}$, define
\[
    R=\left\lceil H/\tileSize\right\rceil,
    \qquad
    C=\left\lceil W/\tileSize\right\rceil .
\]
The operator $\mathcal{B}_{\tileSize}$ pads $\mA$ with zeros, if needed, to an
$(R\tileSize)\times(C\tileSize)$ matrix and partitions the padded matrix into
an $R\times C$ grid of $\tileSize\times\tileSize$ tiles:
\[
    \mathcal{B}_{\tileSize}(\mA)
    =
    \{\mA_{rc}\}_{r=1,\ldots,R;\,c=1,\ldots,C},
    \qquad
    \mA_{rc}\in\R^{\tileSize\times\tileSize}.
\]
The inverse operation $\mathcal{B}_{\tileSize}^{-1}$ denotes reassembly of the
tiles followed by removal of the padding.  We reserve the notation
$\tileSize=\infty$ for the no-tiling case, in which the full matrix is processed
as a single object rather than padded to a square tile.  Rectangular tiles
$(T_h,T_w)$ are used only in the tile-shape study of
\cref{fig:tile_shape_tradeoff}; elsewhere $T_h=T_w=\tileSize$.

\subsection{Polar Factor and Newton--Schulz Map}
\label{sec:polar_ns}

We first recall the row-oriented polar factor.  Suppose
$\mA\in\R^{H\times W}$ has full row rank and $H\le W$.  Its row polar
decomposition is
\begin{equation}
\label{eq:polar_decomp}
    \mA = \mU\mP,
\end{equation}
where $\mU\in\R^{H\times W}$ has orthonormal rows,
$\mU\mU^\top=\mI$, and $\mP=(\mA^\top\mA)^{1/2}$ is symmetric positive
semidefinite~\cite{P1124,P1208}.  The factor
$\mU=\polar(\mA)$ is the closest row-orthonormal matrix to $\mA$ in any
unitarily invariant norm. When $\mA$ is rank deficient, the polar factor need not be unique.
In this paper, $\polar(\mA)$ is used only as an exact reference
matrix function; the optimizer update itself is defined by the finite
Newton--Schulz map below.

Muon uses a quintic Newton--Schulz iteration as a finite polar-factor
surrogate.  For a nonzero matrix $\mA\in\R^{H\times W}$ with $H\le W$, set
\[
    \mX_0=\frac{\mA}{\normF{\mA}} .
\]
The row-oriented recurrence is
\begin{equation}
\label{eq:ns_iteration}
    \mX_{k+1}
    =
    \bigl(
        \nsA\mI
        + \nsB\mX_k\mX_k^\top
        + \nsC(\mX_k\mX_k^\top)^2
    \bigr)\mX_k ,
\end{equation}
with coefficients
\begin{equation}
\label{eq:ns_quintic}
    (\nsA,\nsB,\nsC)=(3.4445,-4.7750,2.0315),
\end{equation}
as in the Muon implementation~\cite{KellerJordan2024muon}.  Here $\mI$ has
dimension $H\times H$.  If $H>W$, the same recurrence is applied to
$\mA^\top$ and the result is transposed back, so the Gram matrix in the
iteration always has dimension $\min\{H,W\}$.  For the zero matrix, we set
$\Phi_K(\mathbf{0})=\mathbf{0}$.

We denote the result after $K$ Newton--Schulz steps by
\begin{equation}
\label{eq:phi_k_def}
    \Phi_K(\mA) := \mX_K .
\end{equation}
The map $\Phi_K$ should not be identified with the exact polar factor.  It is a
finite Newton--Schulz map used to construct the optimizer direction.  Although
it drives singular values toward a common scale, it is not an exact polar
decomposition.

\subsection{The Muon Update}
\label{sec:muon}

We now state the full-matrix Muon update that HiMuon modifies.  Let
$\mW_{t-1}\in\R^{H\times W}$ be a two-dimensional parameter matrix at optimizer
step $t$, and let
\[
    \mG_t = \nabla_{\mW}\Loss(\vtheta_{t-1})
\]
be its gradient.  Muon first updates a momentum buffer,
\begin{equation}
\label{eq:muon_momentum}
    \mM_t = \beta\mM_{t-1}+\mG_t ,
\end{equation}
and then forms the Nesterov-type momentum gradient
\begin{equation}
\label{eq:muon_lookahead}
    \widehat{\mG}_t = \mG_t+\beta\mM_t .
\end{equation}
The matrix-valued update direction is obtained by applying the finite
Newton--Schulz map to this momentum gradient:
\begin{equation}
\label{eq:muon_direction}
    \mU_t = \Phi_K(\widehat{\mG}_t).
\end{equation}
With learning rate $\lr$ and decoupled weight-decay coefficient $\lambda$, the
core parameter update is
\begin{equation}
\label{eq:muon_update}
    \mW_t
    =
    (1-\lr\lambda)\mW_{t-1}
    -
    \lr \mU_t .
\end{equation}
Equations~\eqref{eq:muon_direction}--\eqref{eq:muon_update} define the
full-matrix baseline.  HiMuon keeps the momentum and weight-decay structure of
Muon, but replaces the single full-matrix application of $\Phi_K$ in
\eqref{eq:muon_direction} by tilewise applications of the same finite map.

\section{Related Work}
\label{sec:related}

This section places HiMuon in context with prior work on Muon-type
optimization and efficient matrix-function updates.  We organize the discussion
around five themes: update scaling, structured decompositions of the
Newton--Schulz computation, alternative reductions of the Muon matrix function,
batched GPU execution, and theoretical or diagnostic benchmarks.

\subsection{Scaling and normalization in Muon-type optimizers}

Muon-type optimizers apply a finite Newton--Schulz map to momentum-gradient
matrices, but their practical behavior also depends on how the resulting update
is scaled.  Liu et al.~\cite{P1066} identify weight growth and
shape-dependent update magnitudes as obstacles to scaling Muon, and address
them through decoupled weight decay and a shape-dependent learning-rate scaler.
Other variants modify the normalization before or after the Newton--Schulz map.
MuonEq~\cite{chang2026muoneq} equilibrates rows and columns before applying the
map; Muon+~\cite{zhang2026muonplus} and NorMuon~\cite{P1069} normalize the
update after the map; RMNP~\cite{deng2026rmnp} replaces the map by row-wise
normalization; and TrasMuon~\cite{TrasMuon2026} adds feature-wise trust-region
clipping.  These methods show that Muon-type optimization depends not only on
singular-value homogenization, but also on how update magnitudes are balanced
across rows, neurons, and parameter shapes.

\subsection{Structured decompositions of Newton--Schulz updates}

A closely related line of work reduces the cost of Muon by replacing one
global Newton--Schulz update with smaller structured subproblems.
GLM-5~\cite{GLM5Team2026} applies Muon independently to each attention-head
submatrix, while Group Muon~\cite{zhang2026groupmuon} treats the head-group size
as a tunable hyperparameter; Boreiko et al.~\cite{P1111} study row-wise,
column-wise, and square block splitting whose block sizes follow the tensor-parallel degree.  MuonBP~\cite{P1082} uses
block orthogonalization together with periodic global correction, while
DASH~\cite{modoranu2026dash} demonstrates the efficiency of batched block
updates in Shampoo-type methods.  TEON~\cite{zhang2026teon} takes a different
direction by coupling multiple layers through tensorized orthogonalization.
These methods decompose matrix-valued update computations using units motivated
by model semantics, parallel layouts, correction schemes, or cross-layer tensor
structure.  HiMuon instead uses fixed two-dimensional intra-matrix tiles
whose size is not determined by the sharding layout, making the tile size both
a numerical parameter and a GPU execution parameter.

\subsection{Alternative reductions of the Muon matrix function}

Other approaches reduce Muon cost by changing the matrix being transformed, the
polynomial map, or the matrix-function surrogate.  Dion2~\cite{Ahn2025dion2}
orthogonalizes selected rows with error feedback, while Magma~\cite{Joo2026magma}
uses stochastic masking based on momentum-gradient alignment.
Turbo-Muon~\cite{P1093} reduces the number of Newton--Schulz steps;
PRISM~\cite{YangPRISM2026} adapts the polynomial coefficients to the observed
spectrum; and IFNSO~\cite{HuIFNSO2026} replaces the iteration by a learned
polynomial operator.  MUD~\cite{SouthworthMUD2026} uses triangular whitening
instead of Newton--Schulz, and PolarGrad~\cite{PolarGrad2026} studies
orthogonalized updates from a broader matrix-optimization perspective.  This
line of work modifies the operator used to approximate or replace the
polar-factor surrogate.  HiMuon keeps the finite Newton--Schulz map itself
unchanged, but changes the domain on which it is applied.

\subsection{Batched GPU execution}

Tiled Newton--Schulz computation exposes many independent small dense matrix
products.  This connects HiMuon to batched small-matrix linear algebra, where
grouping independent operations is important for GPU utilization~\cite{P1226}.
Related implementation ideas include epilogue fusion,
on-chip data reuse, and custom Triton kernels~\cite{P1165,P1227}.  HiMuon uses
these ideas inside the optimizer step through tilewise Newton--Schulz calls,
cross-layer batching, memory-bounded chunking, and runtime tile-size
reconfiguration.

\subsection{Theory and diagnostic benchmarks}

Recent theory has begun to clarify why Muon-type updates can differ from
scalar adaptive methods.  Kim and Oh~\cite{KimOhMuonNS2026} analyze
finite-step Muon convergence; Su~\cite{Su2025isotropic} studies
singular-value homogenization under curvature models; Du and
Su~\cite{DuNewtonMuon2026} derive an activation-aware Newton-Muon update; and
Wen et al.~\cite{P1136} benchmark matrix-based optimizers across model scales.
The associative-memory benchmark of Kim et al.~\cite{KimNichaniCapacity2026}
provides a controlled setting in which full-matrix Muon and SGD have sharply
separated one-step capacity scalings.  We use this benchmark in
\Cref{sec:capacity} as a diagnostic for whether tilewise Newton--Schulz updates
retain the qualitative signal-amplification behavior of full-matrix Muon.

Taken together, these lines of work motivate the design studied in this paper.
HiMuon keeps the finite Newton--Schulz map used by Muon, but applies it to
fixed-size two-dimensional intra-matrix tiles.  This choice makes the tile size
a numerical parameter, controlling the amount of retained spectral coupling,
and an execution parameter, producing uniform small dense matrix problems for
batched GPU kernels.  The resulting optimizer combines local matrix-function
updates, runtime tile-size scheduling, and cross-layer batched execution.  The
capacity diagnostic in \Cref{sec:capacity} is then used to assess how this
local replacement affects the qualitative behavior of full-matrix Muon.
\section{The HiMuon Method}
\label{sec:method}

We now present HiMuon.  Starting from the full-matrix Muon
update in \eqref{eq:muon_direction}, HiMuon replaces the single application of
$\Phi_K$ to a momentum-gradient matrix by independent applications of the same
finite Newton--Schulz map on fixed-size tiles.  The tiled update is
then implemented as a batched small-matrix workload on the GPU.

\subsection{Tiled Newton--Schulz Map}
\label{sec:tiling}

\begin{figure}[t]
    \centering
    \includegraphics[width=\textwidth]{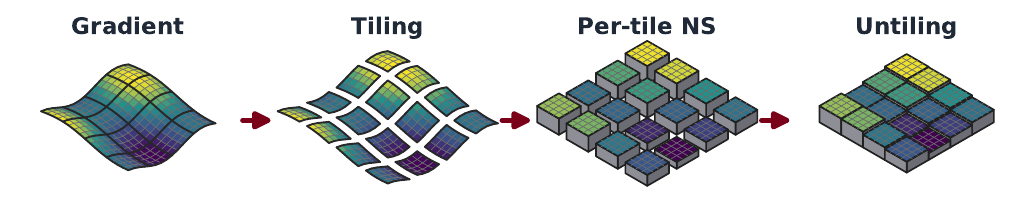}
    \caption{The tiled Newton--Schulz map $\Phi_{K,\tileSize}$.  The lookahead
    momentum gradient $\widehat{\mG}$ is partitioned into
    $\tileSize\times\tileSize$ tiles.  The finite Newton--Schulz map $\Phi_K$
    is applied independently to each tile, and the tile outputs are reassembled
    into $\Phi_{K,\tileSize}(\widehat{\mG})$.  Because each tile is processed in
    isolation, spectral interactions are retained within tiles but removed
    across tile boundaries.}
    \label{fig:tiled_ns}
\end{figure}

Let $\widehat{\mG}\in\R^{H\times W}$ be the lookahead momentum gradient defined
in \eqref{eq:muon_lookahead}.  In full-matrix Muon, the update direction is
$\Phi_K(\widehat{\mG})$, where $\Phi_K$ is the finite Newton--Schulz map defined
in \Cref{sec:polar_ns}.  For a finite tile size $\tileSize$, HiMuon applies this
map locally.  Using the tiling operator $\mathcal{B}_{\tileSize}$ from
\Cref{sec:notation}, we write
\begin{equation}
\label{eq:tiling_operator}
    \mathcal{B}_{\tileSize}(\widehat{\mG})
    =
    \bigl\{\widehat{\mG}_{rc}\bigr\}_{r=1,\ldots,R;\,c=1,\ldots,C},
    \qquad
    R=\left\lceil H/\tileSize\right\rceil,\quad
    C=\left\lceil W/\tileSize\right\rceil ,
\end{equation}
where zero padding is used when needed and each tile
$\widehat{\mG}_{rc}$ has size $\tileSize\times\tileSize$.  HiMuon applies the
same map $\Phi_K$ independently to each tile and then reassembles the outputs:
\begin{equation}
\label{eq:tiled_ns}
    \Phi_{K,\tileSize}(\widehat{\mG})
    :=
    \mathcal{B}_{\tileSize}^{-1}
    \!\left(
        \bigl\{\Phi_K(\widehat{\mG}_{rc})\bigr\}_{r=1,\ldots,R;\,c=1,\ldots,C}
    \right).
\end{equation}
The inverse operator removes any padding introduced during tiling.
\Cref{fig:tiled_ns} illustrates the construction.  The special
case $\tileSize=\infty$ denotes the full-matrix branch,
\begin{equation}
\label{eq:full_branch}
    \Phi_{K,\infty}(\widehat{\mG})
    :=
    \Phi_K(\widehat{\mG}) .
\end{equation}

For finite $\tileSize$ below the matrix dimensions, the tiled map is generally
different from the full-matrix Muon map:
\begin{equation}
\label{eq:tiled_not_global}
    \Phi_{K,\tileSize}(\widehat{\mG})
    \neq
    \Phi_K(\widehat{\mG}) .
\end{equation}
This difference is structural: tiling does not commute with either the exact
polar factor or the finite Newton--Schulz surrogate.  It is a local matrix-function
map whose tile size controls the amount of retained spectral coupling.  Larger
tiles preserve more of the full-matrix interaction, while smaller tiles reduce
the cost and expose more independent work.

For square tiles, each tile requires $\mathcal{O}(\tileSize^3K)$ arithmetic,
and there are approximately $HW/\tileSize^2$ tiles.  Thus, up to lower-order
padding effects, the leading Newton--Schulz work is
\begin{equation}
\label{eq:tiled_cost}
    \mathcal{O}(HW\tileSize K),
\end{equation}
compared with the full-matrix cost
$\mathcal{O}(HW\min\{H,W\}K)$.  The same parameter $\tileSize$ determines both
the numerical locality of the update and the size of the dense matrix problems
issued to the GPU.  

\subsection{Update Scaling}
\label{sec:lr_adjust}
As in Muon, we multiply the Newton--Schulz direction by a shape-dependent
RMS-matching factor before applying the parameter update.  For a full
$H\times W$ matrix, this factor is
$c\sqrt{\max\{H,W\}}$, since a semi-orthogonal reference direction has
entrywise RMS $1/\sqrt{\max\{H,W\}}$.  In the tiled case, the Newton--Schulz
map acts on $\tileSize\times\tileSize$ blocks, so the corresponding factor is
the single tile-size-dependent scalar
\[
    \tau(\tileSize)=c\sqrt{\tileSize}.
\]
For the no-tiling branch, we set
$\tau_\ell(\infty)=c\sqrt{\max\{H_\ell,W_\ell\}}$ for
$\mW^{(\ell)}\in\R^{H_\ell\times W_\ell}$.  The scaled update is therefore
\begin{equation}
\label{eq:scaled_update}
    \mW_t^{(\ell)}
    =
    (1-\lr\lambda)\mW_{t-1}^{(\ell)}
    -
    \lr\,\tau_\ell(\tileSize_t)\,
    \Phi_{K,\tileSize_t}(\widehat{\mG}_t^{(\ell)}),
\end{equation}
where
\[
    \tau_\ell(\tileSize_t)
    =
    \begin{cases}
        c\sqrt{\max\{H_\ell,W_\ell\}}, & \tileSize_t\ge\max\{H_\ell,W_\ell\},\\
        c\sqrt{\tileSize_t}, & \tileSize_t<\max\{H_\ell,W_\ell\} .
    \end{cases}
\]
This scaling is a global scalar for each matrix update; it introduces no
tilewise adaptation and leaves the learning rate, momentum, and weight decay
unchanged from the Muon baseline.

\subsection{Complete Algorithm}
\label{sec:alg1}

\Cref{alg:himuon} summarizes the HiMuon update.  The index
$\ell\in\mathcal{E}$ enumerates the two-dimensional parameter matrices to which
Muon-type updates are applied, with
$\mW^{(\ell)}\in\R^{H_\ell\times W_\ell}$.  The value $\tileSize_t$ denotes the
tile size used at optimizer step $t$; the special value $\tileSize_t=\infty$
selects the full-matrix branch.

\begin{algorithm}[ht]
\caption{HiMuon}
\label{alg:himuon}
\begin{algorithmic}[1]
\REQUIRE learning rate $\lr$, momentum $\beta$, Nesterov flag, weight decay $\lambda$, NS steps $K$, tile-size schedule $\{\tileSize_t\}$, scaling constant $c$
\STATE Initialize $\mM_0^{(\ell)}\gets\mathbf{0}$ for all $\ell\in\mathcal{E}$
\FOR{$t=1,2,\ldots$}
    \STATE Compute gradients
    $\mG_t^{(\ell)} \gets \nabla_{\mW^{(\ell)}}\Loss(\vtheta_{t-1})$
    for all $\ell\in\mathcal{E}$
    \FOR{each $\ell\in\mathcal{E}$}
        \STATE $\mM_t^{(\ell)} \gets
        \beta\mM_{t-1}^{(\ell)} + \mG_t^{(\ell)}$
        \IF{Nesterov}
            \STATE $\widehat{\mG}_t^{(\ell)} \gets
            \mG_t^{(\ell)} + \beta\mM_t^{(\ell)}$
        \ELSE
            \STATE $\widehat{\mG}_t^{(\ell)} \gets
            \mM_t^{(\ell)}$
        \ENDIF

        \IF{$\tileSize_t \ge \max\{H_\ell,W_\ell\}$}
            \STATE $\mU_t^{(\ell)} \gets
            \Phi_K(\widehat{\mG}_t^{(\ell)})$
            \STATE $\tau_t^{(\ell)} \gets
            c\sqrt{\max\{H_\ell,W_\ell\}}$
        \ELSE
            \STATE $\{\widehat{\mG}_{t,rc}^{(\ell)}\}_{r,c}
            \gets
            \mathcal{B}_{\tileSize_t}
            (\widehat{\mG}_t^{(\ell)})$
            \hfill \COMMENT{tile}
            \STATE $\{\mU_{t,rc}^{(\ell)}\}_{r,c}
            \gets
            \{\Phi_K(\widehat{\mG}_{t,rc}^{(\ell)})\}_{r,c}$
            \hfill \COMMENT{tilewise NS}
            \STATE $\mU_t^{(\ell)} \gets
            \mathcal{B}_{\tileSize_t}^{-1}
            \bigl(\{\mU_{t,rc}^{(\ell)}\}_{r,c}\bigr)$
            \hfill \COMMENT{untile}
            \STATE $\tau_t^{(\ell)} \gets c\sqrt{\tileSize_t}$
        \ENDIF

        \STATE $\mW_t^{(\ell)} \gets
        (1-\lr\lambda)\mW_{t-1}^{(\ell)}
        - \lr\,\tau_t^{(\ell)}\mU_t^{(\ell)}$
    \ENDFOR
\ENDFOR
\end{algorithmic}
\end{algorithm}

The branch $\tileSize_t\ge\max\{H_\ell,W_\ell\}$, for example
$\tileSize_t=\infty$, denotes the full-matrix Muon update.  Smaller tile sizes
are interpreted through the padding-and-tiling operator
$\mathcal{B}_{\tileSize_t}$.  For square finite tiles,
$\tau_t^{(\ell)}=c\sqrt{\tileSize_t}$ is a single scalar shared by all tiles in
the matrix; it introduces no tilewise adaptation.  Rectangular tiles use
$\tau_t^{(\ell)}=c\sqrt{\max\{T_h,T_w\}}$ in the tile-shape experiments.

\subsection{GPU Implementation}
\label{sec:hardware}

The mathematical definition of HiMuon creates many independent
Newton--Schulz problems, one for each tile.  This independence is useful only
if the implementation can execute the tilewise computations without introducing
new overheads that offset the arithmetic reduction.  A naive implementation
would allocate a fresh tile tensor and launch many small Newton--Schulz
computations every step, writing intermediate matrices to high-bandwidth memory
(HBM) between the $K$ Newton--Schulz steps.  These costs can dominate when the
tile size is small.

The implementation contribution is therefore an optimizer-level execution plan
for the tiled Newton--Schulz map.  The plan has three goals: avoid unnecessary
data movement, aggregate small independent tile computations into large GPU
batches, and allow the tile size to change during training without resetting
optimizer state.  These choices do not change the mathematical update in
\Cref{alg:himuon}; they determine how the independent evaluations of
$\Phi_K$ are scheduled on the GPU.

\subsubsection{Tile-size-dependent Newton--Schulz kernels}
\label{sec:kernel_paths}

The best kernel strategy depends on the tile size.  For larger tiles, each
Newton--Schulz step contains enough arithmetic to amortize kernel launches, and
HiMuon uses a multi-kernel Triton~\cite{P1165} path.  In this path, matrix
products and the polynomial update in the Newton--Schulz recurrence are fused
where possible before intermediate results are written back to HBM.

For smaller tiles, the arithmetic per tile is lower and repeated memory traffic
becomes the main bottleneck.  In this regime, the tile working set fits in
on-chip static random-access memory (SRAM).  HiMuon therefore uses an
SRAM-resident kernel that keeps all $K$ Newton--Schulz iterations inside a
single kernel, avoiding intermediate HBM round trips through the iteration
loop.  The threshold between the two paths is hardware dependent; in our
implementation, the SRAM-resident path is used for tiles up to
$\tileSize=128$.

\begin{figure}[htbp]
    \centering
    \includegraphics[width=0.85\textwidth]{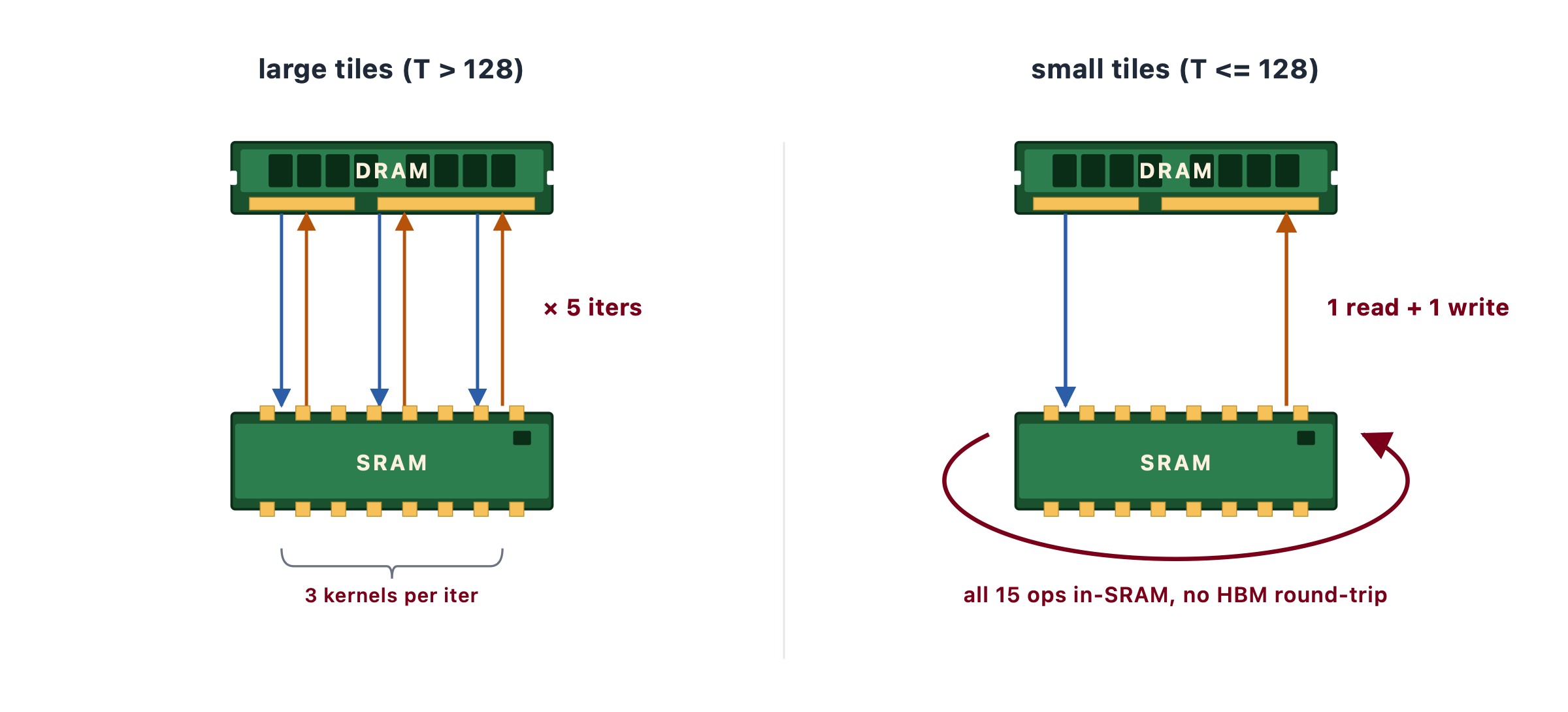}
    \caption{Tile-size-dependent kernel strategy.  Larger tiles use a
    multi-kernel path with epilogue fusion.  Smaller tiles use an
    SRAM-resident path in which all $K$ Newton--Schulz iterations are kept
    on chip, avoiding intermediate HBM round trips.}
    \label{fig:kernel_paths}
\end{figure}

\Cref{fig:kernel_paths} summarizes this tile-size-dependent dispatch rule.
These kernel choices do not change the arithmetic complexity
$\mathcal{O}(HW\tileSize K)$; their purpose is to make the reduced arithmetic
cost visible in wall-clock time by reducing memory traffic and kernel-launch
overhead.

\subsubsection{Cross-layer batching}
\label{sec:xlayer}

Tilewise Newton--Schulz evaluations are independent not only within a single
matrix, but also across different parameter matrices.  This is important in
transformer models, where each layer contains several Muon-eligible matrices
and each matrix may contain only a moderate number of tiles.  Launching one
Newton--Schulz computation per matrix can underutilize the GPU because each
launch may expose too little work.

HiMuon thus groups compatible tiles across layers into larger batched
Newton--Schulz calls.  Tiles are compatible when they have the same tile shape,
dtype, and kernel path.  The batched call computes the same tilewise outputs as
processing each matrix separately; only the execution order changes.  Momentum
buffers, parameter tensors, and weight decay remain indexed by the original
parameter $\ell$, so cross-layer batching is an implementation strategy rather
than a change to the optimizer update.

For large models, a single global batch of all compatible tiles may exceed the
available memory budget.  HiMuon handles this with memory-bounded chunking: the
batched tile set is split into the fewest chunks whose input, intermediate
workspace, and output fit within a prescribed memory budget.  Since the
Newton--Schulz evaluations are independent across tiles, chunking is
mathematically equivalent to processing the full batch at once.  It only bounds
temporary memory usage.

In single-GPU and Distributed Data Parallel (DDP) runs, the batched tiles are assembled during each optimizer step into reusable batched buffers
whose storage is cached across steps.  Under the bank-sharded
Fully Sharded Data Parallel (FSDP) layout used in our experiments, same-shape
matrices are stored along a leading batch axis, so each rank owns a subset of complete
matrices in its local shard.  The tiled Newton--Schulz computation is therefore
rank-local and does not require an optimizer-internal all-gather.

\subsubsection{Runtime tile-size reconfiguration}
\label{sec:schedule}

The tile size is exposed as a runtime parameter rather than being fixed at
optimizer construction.  A \texttt{reconfigure()} routine updates the
tile-dependent execution plan when $\tileSize_t$ changes.  This plan includes
tile indexing metadata, batching buckets, chunking decisions, kernel-path
selection, and graph-capture state when applicable.  The optimizer state itself,
including the momentum buffers, is independent of the tile size and is
preserved across reconfiguration.

We treat $\tileSize_t=\infty$ as the full-matrix branch, in which $\Phi_K$ is
applied to the original matrix.  Finite tile sizes are interpreted through the
padded tiling operator $\mathcal{B}_{\tileSize_t}$, where zero padding only
aligns the matrix boundaries with the tile grid: padded rows and columns
remain identically zero throughout the iteration \eqref{eq:ns_iteration} and
do not affect the Frobenius normalization of the tile.  Tiling is therefore
meaningful only when $\tileSize_t<\max\{H_\ell,W_\ell\}$.  When
$\tileSize_t\ge\max\{H_\ell,W_\ell\}$, the partition reduces to a single tile
($R=C=1$) and the tiled map coincides with the full-matrix map,
$\Phi_{K,\tileSize_t}(\widehat{\mG}_t^{(\ell)})
=\Phi_K(\widehat{\mG}_t^{(\ell)})$; this case is therefore identified with the
full-matrix branch, including its scaling factor
$c\sqrt{\max\{H_\ell,W_\ell\}}$.

In the end-to-end experiments of \Cref{sec:exp_cross_arch}, we use a
full-to-local schedule,
\[
    \tileSize_t
    =
    \begin{cases}
        \infty, & t < t_1,\\
        512,    & t_1 \le t < t_2,\\
        128,    & t \ge t_2 .
    \end{cases}
\]
The initial full-matrix phase preserves global spectral coupling during the
earliest optimization steps, when gradient statistics and dominant directions
are least settled.  The optimizer then switches to lower-cost tiled maps for
the remainder of training.  The DDP experiments use
$(t_1,t_2)=(100,500)$, while the FSDP experiment uses
$(t_1,t_2)=(200,500)$.  These breakpoints are fixed empirical choices in our
experiments.  The capacity diagnostic in \Cref{sec:capacity} motivates the
full-to-local form of the schedule qualitatively, but does not determine the
switch times.
\section{Capacity Diagnostic}
\label{sec:capacity}

Many-step transformer training entangles the effect of the update map with
learning-rate schedules, changing gradient statistics, stochasticity, and
architecture-specific dynamics.  To isolate the effect of replacing a
full-matrix Newton--Schulz map by a tiled one, we use the one-step linear
associative-memory benchmark of Kim et al.~\cite{KimNichaniCapacity2026}.  This
benchmark is useful because full-matrix Muon and SGD have sharply separated
capacity scalings, so the signal-amplification behavior of matrix-function
updates can be observed directly.  We use the benchmark as an empirical
diagnostic for how the tiled map $\Phi_{K,\tileSize}$ changes this behavior as
the tile size varies.

\subsection{Linear Associative-Memory Benchmark}
\label{sec:capacity_setup}

We first recall the benchmark.  The task consists of $N$ key--value pairs
$(v_i,u_i)$, where
\[
    v_i,u_i \sim \mathcal{N}(0,d^{-1}\mI_d).
\]
The associations are stored in a matrix $\mW\in\R^{d\times d}$.  Given an
input key $v_i$, the model assigns logits
\[
    u_j^{\!\top}\mW v_i,\qquad j=1,\ldots,N,
\]
and is trained with cross-entropy loss to predict the matching value index
$i$.  The $i$th training example is sampled with probability
$p_i\propto i^{-\alpha}$.  Starting from $\mW_0=0$, one optimizer update using a
minibatch of size $B$ produces $\mW_1$.  An item $k$ is counted as recovered if
\[
    \arg\max_j u_j^{\!\top}\mW_1 v_k = k .
\]
The one-step storage capacity is the number of recovered items.

The benchmark separates scalar gradient updates from spectral matrix-function
updates through their capacity scalings.  Under the assumptions of
Kim et al.~\cite{KimNichaniCapacity2026}, with $\alpha>1$ and appropriate
scaling of the stabilized polar map, the full-matrix Muon update has one-step
capacity
\begin{equation}
\label{eq:kim_muon_capacity}
    \widetilde{\Theta}
    \!\left(
        \min\{d^{1+1/(2\alpha)},\,B^{1/\alpha}\}
    \right),
\end{equation}
where $B$ is the minibatch size.  Under the same data model, SGD has capacity
\begin{equation}
\label{eq:kim_sgd_capacity}
    \widetilde{\Theta}
    \!\left(
        \min\{d^{1/(2\alpha)},\,B^{1/\alpha}\}
    \right).
\end{equation}
Thus, in the dimension-limited regime, the Muon capacity term is larger than
the SGD term by a factor of $d$.  The shared term $B^{1/\alpha}$ is a
finite-sample ceiling: items absent from the minibatch provide no direct
training signal in the one-step gradient.

A useful interpretation of \eqref{eq:kim_muon_capacity} is bulk
singular-value amplification.  At initialization, the one-step gradient
contains strong directions associated with high-frequency sampled items and a
bulk of weaker directions associated with lower-frequency items.  A stabilized
polar-type update amplifies this bulk more uniformly than a scalar gradient
update, enabling recovery of weaker associations.  Since HiMuon applies the
finite Newton--Schulz map only within tiles, the experiments below test how
much of this full-matrix amplification behavior remains after tilewise
processing and reassembly.

\subsection{Tile-Size Capacity Diagnostic}
\label{sec:capacity_diagnostic}

The capacity scalings in \eqref{eq:kim_muon_capacity} and
\eqref{eq:kim_sgd_capacity} are stated for global update maps.  HiMuon replaces
the full map $\Phi_K(\cdot)$ by the tiled map
$\Phi_{K,\tileSize}(\cdot)$.  Each tile can amplify only the singular
directions present in its own submatrix, while singular directions whose
support spans several tiles are not treated coherently by a single
Newton--Schulz map.  This suggests a tile-size heuristic: larger tiles should
better preserve full-matrix behavior because each tile contains a richer local
spectrum, whereas smaller tiles reduce cost but may weaken the reassembled
signal amplification.  We therefore evaluate two one-step diagnostics:
aggregate capacity scaling with dimension and top-$k$ recovery over the
sampled item support.

\begin{figure}[htbp]
    \centering
    \includegraphics[width=\linewidth]{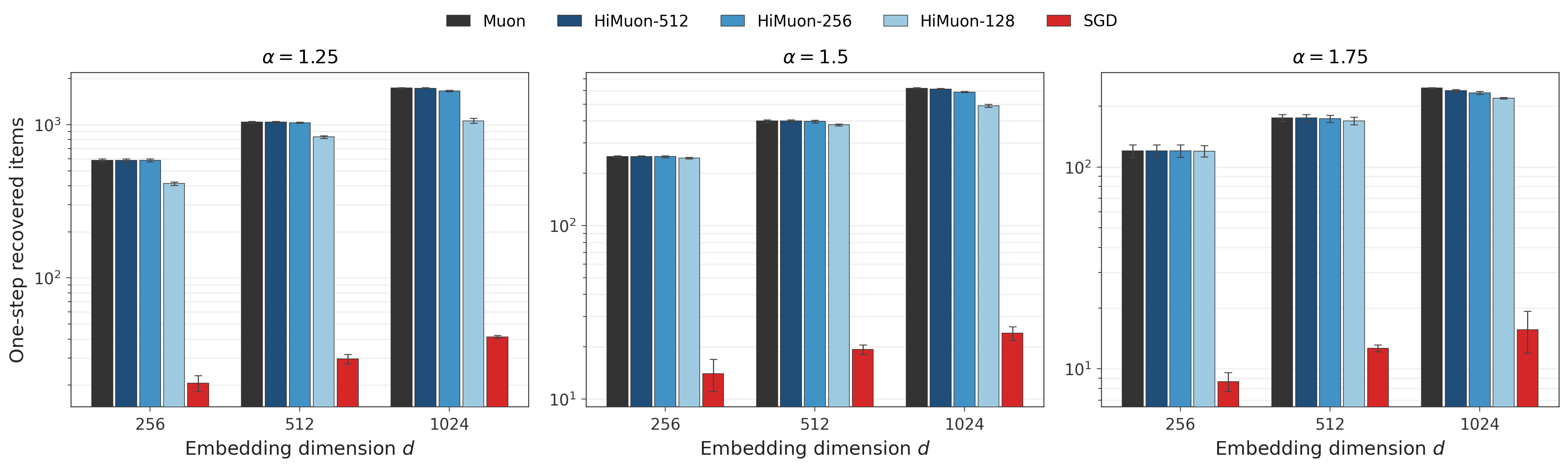}
    \caption{One-step capacity on the associative-memory benchmark at
    $\alpha\in\{1.25,1.5,1.75\}$, $B/d=10$, $N=10^5$, one update from
    $\mW_0=0$, and three random seeds.  Bars are grouped by dimension $d$, and
    error bars show seed standard deviation.  HiMuon approaches full-matrix
    Muon as $\tileSize$ grows; all Newton--Schulz-based updates remain well
    above SGD, and the Muon/SGD gap widens with $d$.}
    \label{fig:capacity_scaling}
\end{figure}

\Cref{fig:capacity_scaling} reports one-step capacity for SGD, full-matrix
Muon, and HiMuon with $\tileSize\in\{128,256,512\}$ across
$\alpha\in\{1.25,1.5,1.75\}$, all within the $\alpha>1$ regime of
\eqref{eq:kim_muon_capacity}.  We use $d\in\{256,512,1024\}$, $B/d=10$,
$N=10^5$, and three random seeds.  Across the tested values of $(d,\alpha)$,
the same qualitative ordering is observed:
\[
    \text{Muon}\approx \text{HiMuon-}512
    >
    \text{HiMuon-}256
    >
    \text{HiMuon-}128
    \gg
    \text{SGD}.
\]
At $d=1024$, HiMuon-$512$ is within $3\%$ of full-matrix Muon across all three
values of $\alpha$, and HiMuon-$256$ is within $6\%$.  HiMuon-$128$ trails
full-matrix Muon by $11\%$--$39\%$, but remains well above SGD in every panel.
The Muon/SGD gap also grows with $d$, consistent with the dimension-dependent
separation predicted by \eqref{eq:kim_muon_capacity} and
\eqref{eq:kim_sgd_capacity}.

\begin{figure}[htbp]
    \centering
    \includegraphics[width=0.875\linewidth]{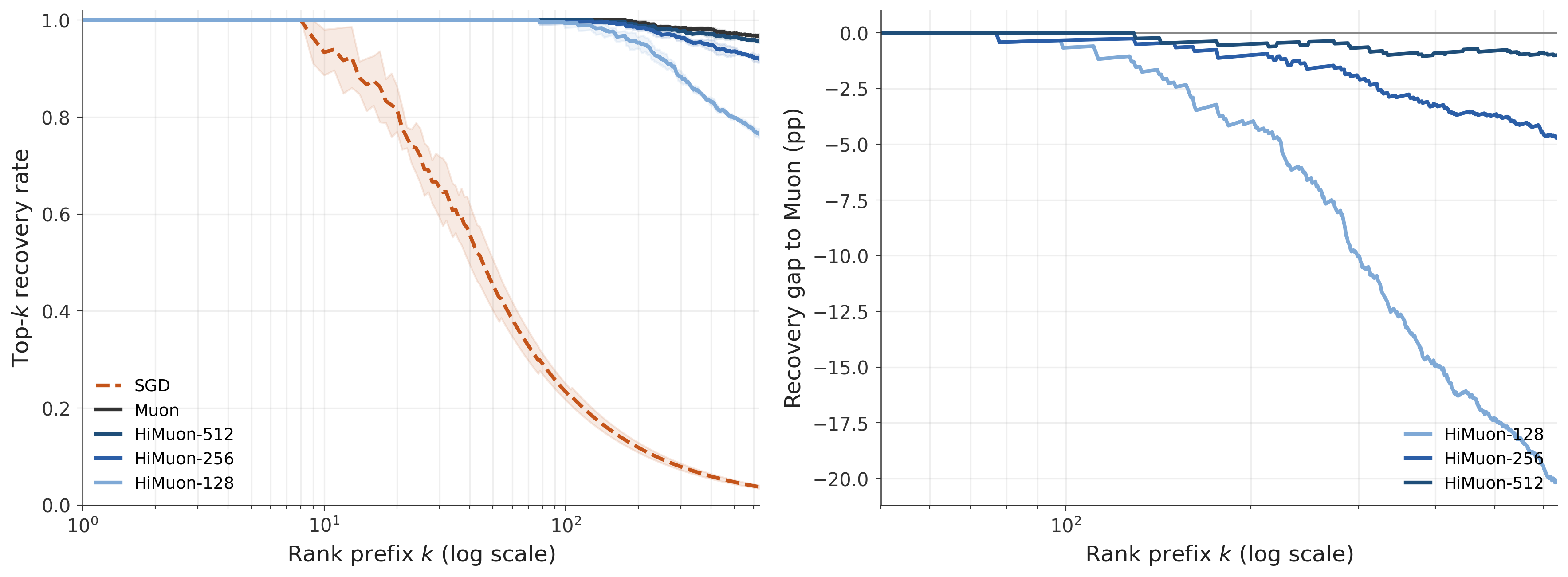}
    \caption{One-step top-$k$ recovery at $d=1024$, $\alpha=1.5$, and
    $B/d=10$.  All methods use the same embeddings, minibatch, and initial
    gradient; only the update map differs.  Recovery is computed over the
    unique sampled support $\mathcal{S}$, sorted by decreasing item frequency.
    \emph{Left}: top-$k$ recovery rate.  \emph{Right}: gap to full-matrix Muon
    in percentage points.  Tile-size differences are small for high-frequency
    items and increase deeper in the sampled support. 
    }
\label{fig:topk_recovery}
\end{figure}

\Cref{fig:topk_recovery} examines which sampled items are recovered at
$d=1024$, $\alpha=1.5$, and $B/d=10$.  We restrict this evaluation to the
unique sampled support $\mathcal{S}$ in the minibatch.  This avoids mixing
update-map quality with finite-sample coverage, since items outside
$\mathcal{S}$ do not contribute directly to the one-step gradient.  Across
seeds, $|\mathcal{S}|$ averages $639$.  We sort sampled items by decreasing
population frequency and report top-$k$ recovery up to the common cutoff
$k_{\max}=632$.

The head of the sampled distribution is insensitive to tile size.  Through
$k\approx128$, full-matrix Muon and all HiMuon variants recover nearly all
sampled top-$k$ items, whereas SGD drops much earlier.  Tile size matters more
deeper in the sampled support, where recovery depends on amplifying weaker
directions.  At $k=k_{\max}$, HiMuon-$128$ trails full-matrix Muon by about
$20$ percentage points, HiMuon-$256$ by about $5$ percentage points, and
HiMuon-$512$ by about $1$ percentage point.  This rank-dependent pattern
supports the interpretation that high-frequency sampled items are robust to
tiling, while recovery of weaker sampled items depends more strongly on the
amount of retained spectral coupling.

These diagnostics suggest that the importance of tile size depends on which
capacity limit is active.  If the finite-sample term $B^{1/\alpha}$ is the
bottleneck, recovery is limited mainly by minibatch coverage, and smaller tiles
may be sufficient.  If the dimension-dependent Muon term in
\eqref{eq:kim_muon_capacity} is the bottleneck, the update map matters more,
and larger tiles are needed to better preserve full-matrix spectral behavior.
In our experiments, $\tileSize=512$ is closest to full-matrix Muon,
$\tileSize=256$ retains much of the gain at lower cost, and
$\tileSize=128$ shows visible capacity loss.

The diagnostic also gives a qualitative rationale for the runtime tile-size
schedule in \Cref{sec:schedule}.  Because the benchmark isolates a one-step
update from $\mW_0=0$, it is most analogous to the earliest stage of training,
where global spectral coupling is most directly exposed.  This motivates using
a full-matrix warm start before switching to finite tiles, and using larger
tiles before smaller ones in a full-to-local schedule.  The benchmark does not
determine the switch times or the final tile size; those choices are fixed
empirical hyperparameters in the end-to-end experiments.
\section{Numerical Experiments}
\label{sec:experiments}

We evaluate HiMuon along four axes: end-to-end training behavior,
matrix-function runtime on transformer layer shapes, optimizer-step speed under
cross-layer batching, and robustness to low-precision Newton--Schulz
computation.  Unless otherwise stated, the baseline is full-matrix
Muon~\cite{P1066}, and HiMuon uses the same learning rate as the corresponding
Muon run.  This isolates the effect of replacing the full-matrix
Newton--Schulz map by the tiled map.

\subsection{Experimental Setup}
\label{sec:exp_setup}

\noindent\textbf{Models and data.}
End-to-end training experiments use Qwen3-0.6B, Qwen3-1.7B, and Qwen3-4B
transformer language models~\cite{P1253}.  The two smaller models are trained
with Distributed Data Parallel (DDP), while the 4B model is trained with Fully
Sharded Data Parallel (FSDP).  Layer-level performance studies use the same
model family to probe the larger-matrix regime.  Training uses the
FineWeb dataset~\cite{P1147} in streaming mode with sequence length $1024$,
drawing from its official \emph{sample-10BT} subset (${\sim}10$B tokens).
Held-out evaluation uses FineWeb's \emph{CC-MAIN-2024-51} configuration, a
Common Crawl snapshot disjoint from the training sample.

\smallskip
\noindent\textbf{Optimizer settings.}
Unless otherwise stated, HiMuon uses momentum $\beta=0.95$ and $K=5$
Newton--Schulz steps, following the default setting of Muon~\cite{P1066}, with
the update scaling of \Cref{sec:lr_adjust} and $c=0.2$.  The end-to-end runs
use the runtime tile-size schedule described in \Cref{sec:schedule}: an initial
full-matrix phase, followed by $\tileSize=512$, then $\tileSize=128$.  The DDP
runs switch at steps $100$ and $500$, while the FSDP run switches at steps
$200$ and $500$.  Isolated matrix-function, kernel, and optimizer-step studies
use a fixed tile size, stated with each experiment.

The learning rate is selected by grid search for full-matrix Muon, and HiMuon
reuses the same value.  The DDP runs use a cosine learning-rate schedule with a
$100$-step linear warmup, decaying to zero.  The FSDP run uses a
warmup-stable-decay (WSD) schedule with the same $100$-step warmup and a final
$400$-step decay to $0.1\times$ the peak learning rate.

\smallskip
\noindent\textbf{Hardware and timing protocol.}
End-to-end training runs use $4\times$ NVIDIA A40 GPUs: Qwen3-0.6B and
Qwen3-1.7B in DDP, and Qwen3-4B in FSDP.  Matrix-function, kernel, and
optimizer-step performance studies run on a single NVIDIA A40 GPU unless
otherwise stated.  Latencies are measured with CUDA events after a warmup
phase, and we report medians with interquartile ranges over $30$--$50$ timed
repetitions.  Layer-level speedups are measured in paired trials: the baseline
and test configuration are run back-to-back, and the per-trial ratio is taken.
Aggregate step-time and kernel speedups are reported as ratios of the
corresponding medians.  Peak memory is the maximum allocated memory
(\texttt{torch.cuda.max\_memory\_allocated}) recorded over a short run of
optimizer steps after warmup.  All experiments use PyTorch~2.6.0 with CUDA~12.4; for details, please refer to our code repository.

\smallskip
\noindent\textbf{Numerical precision.}
All training and performance experiments use bf16 mixed-precision computation,
matching the default Muon implementation.  The precision studies in
\Cref{sec:precision_recovery,sec:precision_descent} vary the precision used
inside the Newton--Schulz inner GEMMs as described there.

\smallskip
\noindent\textbf{Code availability.}
The implementation and the scripts for all experiments in this section are
available at \url{https://github.com/tang0389/himuon}.

\subsection{Training Quality}
\label{sec:exp_cross_arch}

We first evaluate the effect of the tiled update on end-to-end pretraining.
HiMuon replaces one global Newton--Schulz map by many local maps, so the goal
is to test whether the wall-clock reduction is accompanied by a visible loss
penalty in short training runs.  We report next-token cross-entropy on both the
training stream and a held-out validation split disjoint from the training
data.

Each model is trained for $1000$ steps on the FineWeb sample-10BT, comparing HiMuon
against full-matrix Muon at the same learning rate and a single fixed seed.
Qwen3-0.6B and Qwen3-1.7B run in $4\times$A40 DDP with an effective batch of
$64$ length-$1024$ sequences ($31.3$M non-padding tokens), using
$\lr=0.02$ and $\lr=0.01$, respectively.  Qwen3-4B runs in $4\times$A40 FSDP
with an effective batch of $32$ length-$512$ sequences and $\lr=0.005$.  The
HiMuon schedule starts with the full-matrix branch, so HiMuon and Muon share
the same update early in training and differ only after the tile-size switches.

\Cref{tab:train_loss} reports the training loss at selected steps, smoothed by
a trailing mean over the preceding $50$ steps to reduce minibatch noise.  Across
checkpoints, the HiMuon--Muon differences are small relative to the loss scale;
the tiled run is sometimes lower and sometimes higher.

\begin{table}[htbp]
    \centering
    \footnotesize
    \caption{Training loss, reported as a trailing mean over the preceding
    $50$ steps, for Muon and HiMuon.  Values are comparable within each model
    but not across models.  The $0.6$B and $1.7$B models use DDP; the $4$B
    model uses FSDP.}
    \label{tab:train_loss}
    \begin{tabular}{llcccccc}
        \toprule
        & & \multicolumn{6}{c}{Step} \\
        \cmidrule(lr){3-8}
        Model & Optimizer & $100$ & $200$ & $300$ & $500$ & $700$ & $1000$ \\
        \midrule
        \multirow{2}{*}{0.6B}
            & Muon   & $3.220$ & $2.913$ & $2.744$ & $2.521$ & $2.338$ & $2.236$ \\
            & HiMuon & $3.183$ & $2.901$ & $2.723$ & $2.493$ & $2.281$ & $2.162$ \\
        \midrule
        \multirow{2}{*}{1.7B}
            & Muon   & $3.169$ & $2.909$ & $2.746$ & $2.507$ & $2.305$ & $2.189$ \\
            & HiMuon & $3.162$ & $2.914$ & $2.736$ & $2.495$ & $2.277$ & $2.142$ \\
        \midrule
        \multirow{2}{*}{4B}
            & Muon   & $4.838$ & $4.590$ & $4.082$ & $3.701$ & $3.698$ & $3.719$ \\
            & HiMuon & $4.836$ & $4.583$ & $4.103$ & $3.721$ & $3.717$ & $3.733$ \\
        \bottomrule
    \end{tabular}
\end{table}

\begin{figure}[htbp]
    \centering
    \includegraphics[width=\textwidth]{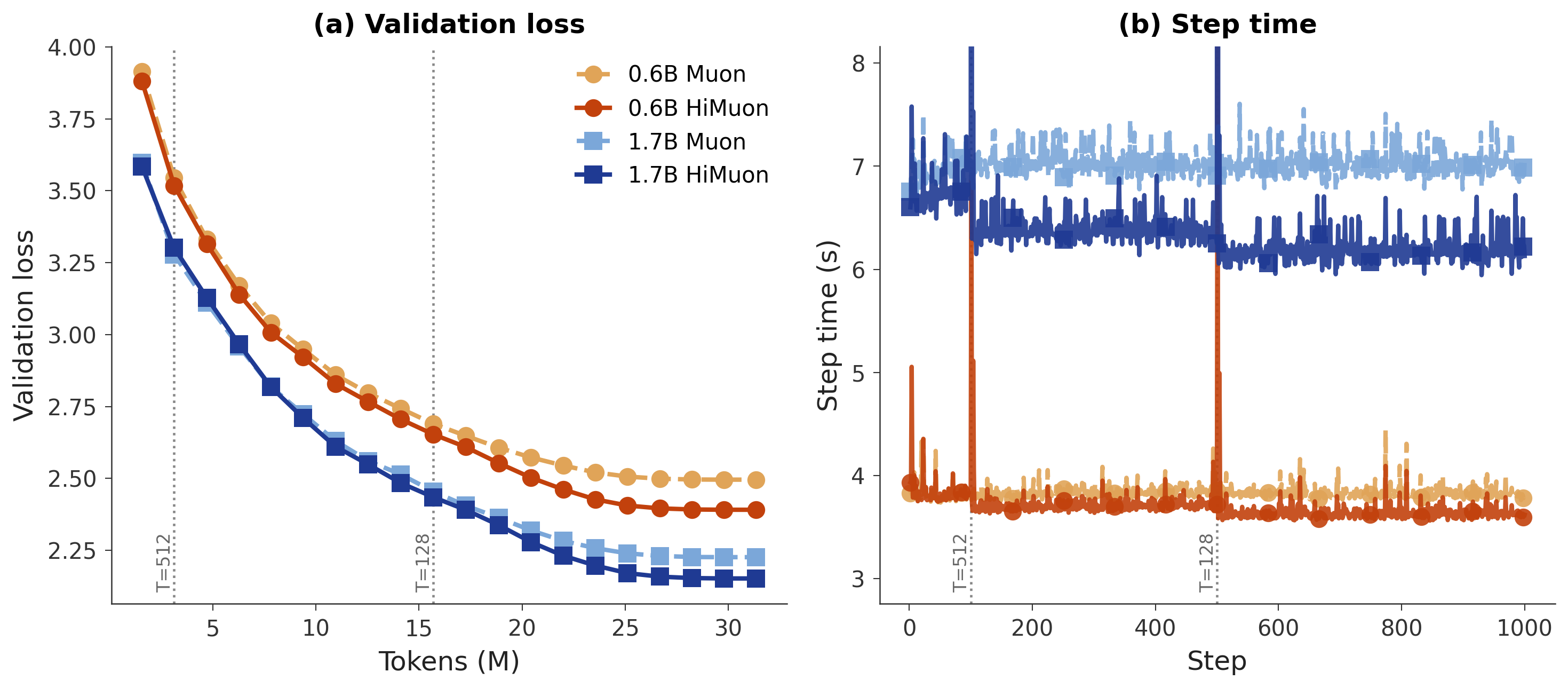}
    \caption{End-to-end DDP training of Qwen3-0.6B and Qwen3-1.7B with HiMuon
    and full-matrix Muon on $4\times$A40 ($1000$ steps, $31.3$M non-padding
    tokens; $\lr=0.02$ and $\lr=0.01$).  \textbf{(a)} Held-out validation loss
    against tokens processed.  \textbf{(b)} Per-step wall time, with the
    tile-size switches marked at step $100$ ($\tileSize=512$) and step $500$
    ($\tileSize=128$).  HiMuon settles to a lower steady-state step time after
    each switch; upward spikes at the switches reflect rebuilding kernels and
    batching metadata.}
    \label{fig:e2e_ddp}
\end{figure}

\smallskip
\noindent\textbf{DDP.}
\Cref{fig:e2e_ddp} plots validation loss against tokens processed and per-step
wall time for the two DDP models.  HiMuon tracks full-matrix Muon through the
initial full-matrix phase and then settles to a lower step time after the
tile-size switches.  In these single-seed runs, the validation losses remain
close to, and in these cases below, the corresponding Muon curves.

For Qwen3-0.6B, the step time decreases from $3.80$\,s in the full-matrix
phase to $3.70$\,s after the first switch and $3.62$\,s after the second switch,
a $4.7\%$ per-step reduction.  Total wall time changes from $3959$\,s to
$3822$\,s, a $3.5\%$ reduction.  The final training loss is lower for HiMuon
by $\Delta=-0.074$ at step $1000$, and the final validation loss is $2.391$
versus $2.496$ for Muon.

For Qwen3-1.7B, the larger matrices make the tiled Newton--Schulz computation
more beneficial.  Step time decreases from $6.72$\,s to $6.36$\,s and then
$6.16$\,s, an $8.3\%$ per-step reduction.  Total wall time changes from
$7157$\,s to $6610$\,s, a $7.6\%$ reduction.  The final training-loss
difference is $\Delta=-0.047$, and the final validation loss is $2.152$ for
HiMuon versus $2.226$ for Muon.  The brief upward spikes in step time coincide
with tile-size switches, where the kernels and cross-layer batching plan are
rebuilt; this overhead is amortized over subsequent steps.

\smallskip
\noindent\textbf{FSDP.}
The Qwen3-4B model is too large for the DDP setup on $4\times$A40 and is
therefore run under FSDP.  Because the Muon baseline in this setting is an
FSDP adaptation rather than the implementation used in the DDP runs, we treat
the 4B experiment as a quality and scaling check rather than a direct timed
comparison.  HiMuon's training loss stays within a few hundredths of
full-matrix Muon in \Cref{tab:train_loss}; the final validation losses are
$3.629$ for HiMuon and $3.613$ for Muon.  HiMuon's step time decreases as the
schedule shrinks the tiles, from $17.55$\,s to $17.32$\,s and then $17.18$\,s
after the switches at steps $200$ and $500$.  End-to-end throughput tuning
under FSDP is left to future work.

Overall, the tiled update gives lower wall-clock time in the DDP runs while
keeping training and validation losses close to full-matrix Muon in these
single-seed experiments.  The reported wall-clock reductions should be
interpreted as short-run measurements: in these $1000$-step runs, the initial
full-matrix phase and schedule-switch overheads form a larger fraction of total
time than they would in longer pretraining runs.

\subsection{Matrix-Function Performance}
\label{sec:exp_matrix_function}

We next study the tiled Newton--Schulz map in isolation, separated from the
full training loop.  This isolates the matrix-function computation from the
forward pass, backward pass, communication, and data movement that also
contribute to end-to-end training time.  The experiments in this subsection
measure speedup on real Qwen3 layer shapes, compare kernel strategies, and
examine the speed--quality trade-off across tile shapes.

\smallskip
\noindent\textbf{Speedup on real transformer layer shapes.}
We first measure the isolated tiled/full-matrix Newton--Schulz speedup on
actual Qwen3 weight matrices.  For each Qwen3 model, we extract the seven
Muon-eligible matrices from Block~0:
\texttt{q\_proj}, \texttt{k\_proj}, \texttt{v\_proj}, \texttt{o\_proj},
\texttt{gate\_proj}, \texttt{up\_proj}, and \texttt{down\_proj}; their shapes
are listed in \Cref{tab:qwen3_shapes}.  For each shape, we apply the
full-matrix map $\Phi_K$ and the tiled map $\Phi_{K,\tileSize}$ with
$\tileSize=512$ and $K=5$ to a single random bf16 Gaussian matrix of that
shape.  We report the paired speedup
$t_{\mathrm{full}}/t_{\mathrm{tile}}$ as the median over $30$ timed repetitions
after $3$ warmup iterations.

\begin{figure}[htbp]
    \centering
    \includegraphics[width=0.75\textwidth]{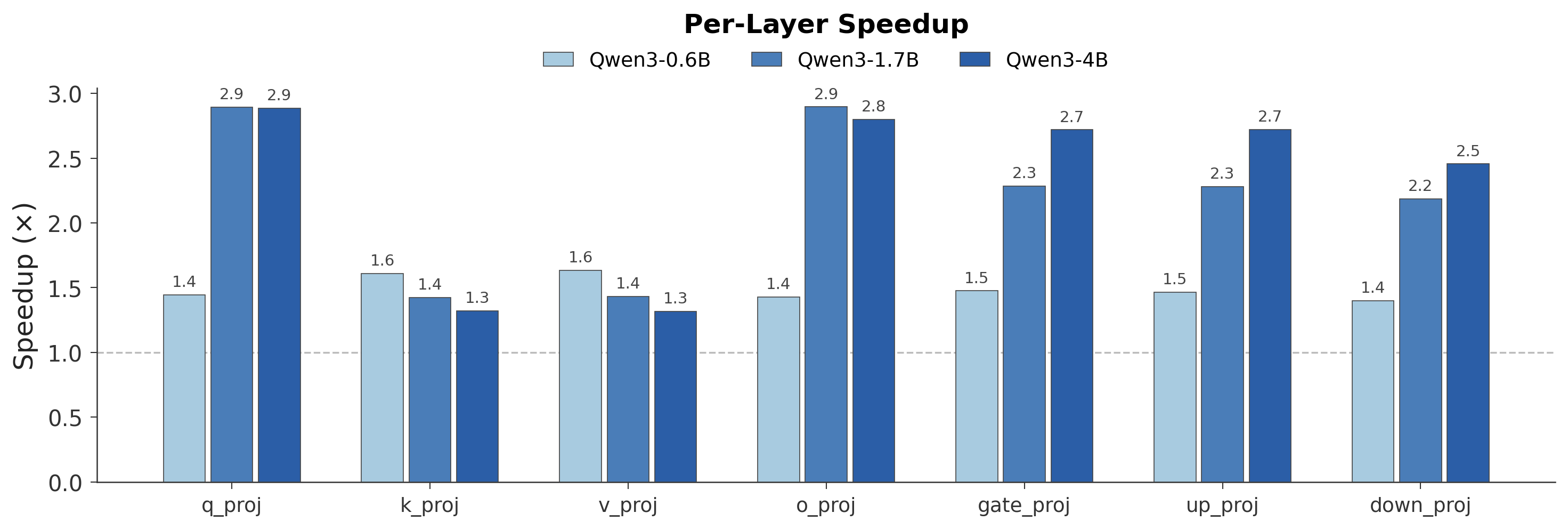}
    \caption{Per-layer tiled/full-matrix Newton--Schulz speedup on real Qwen3
    weight shapes from Block~0 with $\tileSize=512$.  Speedups reach up to
    $2.9\times$ on large attention and MLP projection matrices in Qwen3-1.7B
    and Qwen3-4B.  The \texttt{k\_proj} and \texttt{v\_proj} matrices have
    smaller output dimensions under grouped-query attention and therefore show
    smaller gains.}
    \label{fig:real_shapes}
\end{figure}

\begin{table}[htbp]
    \centering
    \caption{Shapes (output\,$\times$\,input) of the seven Muon-eligible weight
    matrices per transformer block of each Qwen3 model~\cite{P1253}, as used in
    \Cref{fig:real_shapes}.}
    \label{tab:qwen3_shapes}
    \setlength{\tabcolsep}{4pt}
    \resizebox{\textwidth}{!}{%
    \begin{tabular}{lccccccc}
        \toprule
        Model & \texttt{q\_proj} & \texttt{k\_proj} & \texttt{v\_proj} &
        \texttt{o\_proj} & \texttt{gate\_proj} & \texttt{up\_proj} &
        \texttt{down\_proj} \\
        \midrule
        0.6B & $2048\times1024$ & $1024\times1024$ & $1024\times1024$ &
        $1024\times2048$ & $3072\times1024$ & $3072\times1024$ &
        $1024\times3072$ \\
        1.7B & $2048\times2048$ & $1024\times2048$ & $1024\times2048$ &
        $2048\times2048$ & $6144\times2048$ & $6144\times2048$ &
        $2048\times6144$ \\
        4B & $4096\times2560$ & $1024\times2560$ & $1024\times2560$ &
        $2560\times4096$ & $9728\times2560$ & $9728\times2560$ &
        $2560\times9728$ \\
        \bottomrule
    \end{tabular}%
    }
\end{table}

\Cref{fig:real_shapes} shows that the measured per-layer speedups at
$\tileSize=512$ range from $1.4$--$1.6\times$ for Qwen3-0.6B,
$1.4$--$2.9\times$ for Qwen3-1.7B, and $1.3$--$2.9\times$ for Qwen3-4B.  The
largest gains occur on the attention projections \texttt{q\_proj} and
\texttt{o\_proj} and on the wide MLP projections \texttt{gate\_proj},
\texttt{up\_proj}, and \texttt{down\_proj}.  In contrast, \texttt{k\_proj} and
\texttt{v\_proj} show smaller gains.  Qwen3 uses grouped-query attention, which
reduces the output dimension of these two projections relative to the hidden
size; the full-matrix Newton--Schulz work on these narrower matrices is
already smaller, so tiling has less work to reduce.  Across the remaining
projections, the measurements show that the $\mathcal{O}(HW\tileSize K)$ cost
model translates into practical speedups on real transformer weight matrices,
with larger gains when $H$ and $W$ are large relative to $\tileSize$.

\smallskip
\noindent\textbf{Kernel strategy.}
The tiled map changes the performance regime of Newton--Schulz.  At small tile
sizes, the arithmetic work per tile is low, and repeated memory traffic can
dominate.  Without fusion, each Newton--Schulz iteration writes intermediate
matrices to HBM and reads them back.  Small tiles are also the regime in which
the full $K$-step computation can fit in on-chip SRAM, enabling a single
SRAM-resident kernel.

We study this regime at $\tileSize=128$, comparing four implementations: a
PyTorch eager loop, \texttt{torch.compile}, a three-kernel Triton
implementation with epilogue fusion, and an SRAM-resident kernel that keeps
all $K=5$ Newton--Schulz iterations on chip.  Each implementation processes a
batch of random bf16 Gaussian $\tileSize\times\tileSize$ tiles.  We sweep the
tile batch size $B_{\mathrm{tile}}$ and report the median over $50$ timed
repetitions after $10$ warmup iterations, with speedup measured relative to
the eager loop.

\begin{figure}[htbp]
    \centering
    \includegraphics[width=0.85\textwidth]{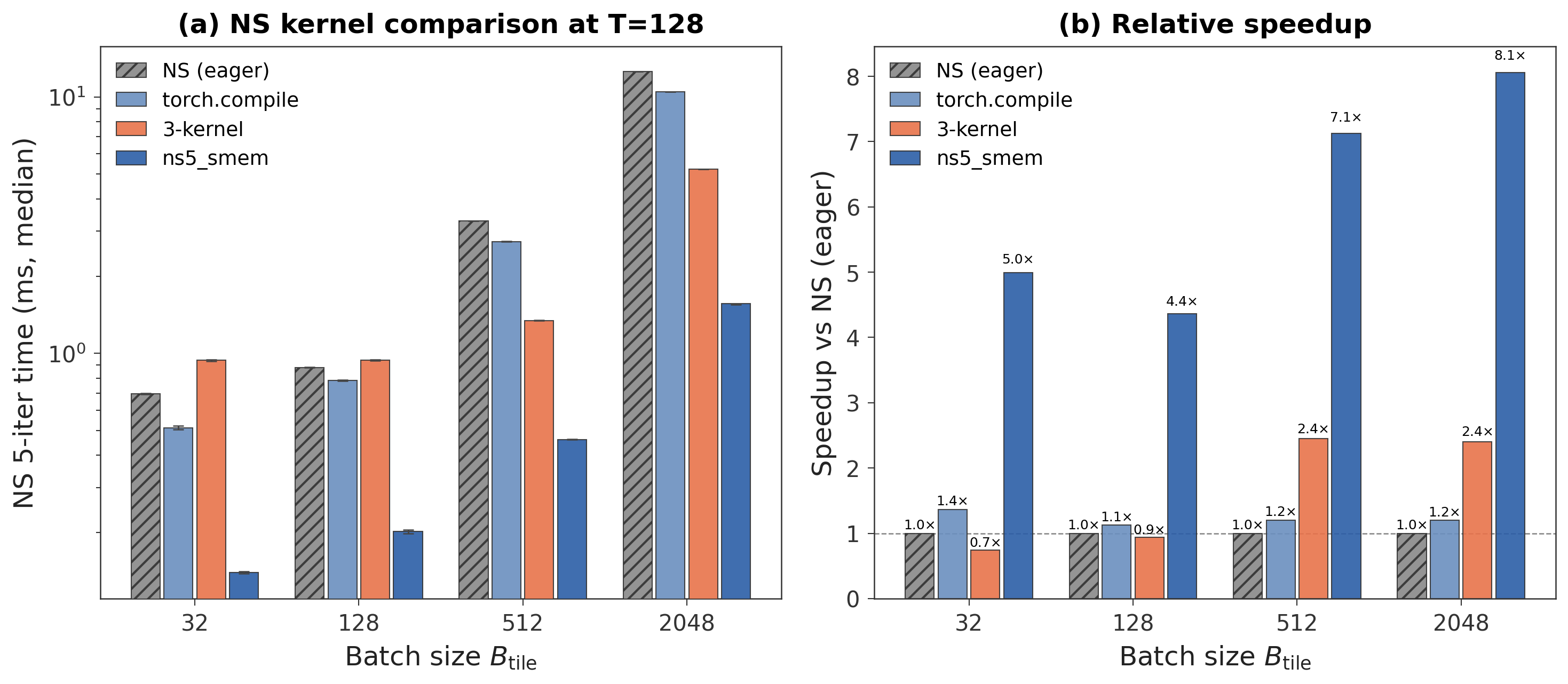}
    \caption{Kernel strategy comparison at $\tileSize=128$ on the NVIDIA A40
    GPU, sweeping tile batch size
    $B_{\mathrm{tile}}\in\{32,128,512,2048\}$.  \textbf{(a)} Absolute
    wall-clock time, shown on a log scale.  \textbf{(b)} Speedup relative to
    the PyTorch eager implementation.  The SRAM-resident kernel removes
    intermediate HBM traffic within the fused NS inner loop and obtains roughly $4$--$8\times$ speedup over
    eager execution across the swept batch sizes.}
    \label{fig:ns_kernel_bench}
\end{figure}

\Cref{fig:ns_kernel_bench} shows that the SRAM-resident path gives the largest
improvement for $\tileSize=128$, achieving roughly $4$--$8\times$ speedup over
eager execution, with the largest gain at
$B_{\mathrm{tile}}=2048$.  The three-kernel epilogue-fusion path is useful
only at larger batch sizes, where its launch overhead is amortized.
\texttt{torch.compile} gives a smaller but consistent improvement, reflecting
partial reduction of launch overhead without fusing the full Newton--Schulz
iteration loop.

\smallskip
\noindent\textbf{Tile shape.}
The tile shape affects both execution and optimization behavior because it
defines the local matrix to which $\Phi_K$ is applied.  We compare rectangular
tile shapes using two measurements: isolated kernel speed and final training
loss for Llama-220M, a small non-Qwen model for inexpensive tile-shape
sweeping, trained with HiMuon for $300$ steps.

\begin{figure}[htbp]
    \centering
    \includegraphics[width=\textwidth]{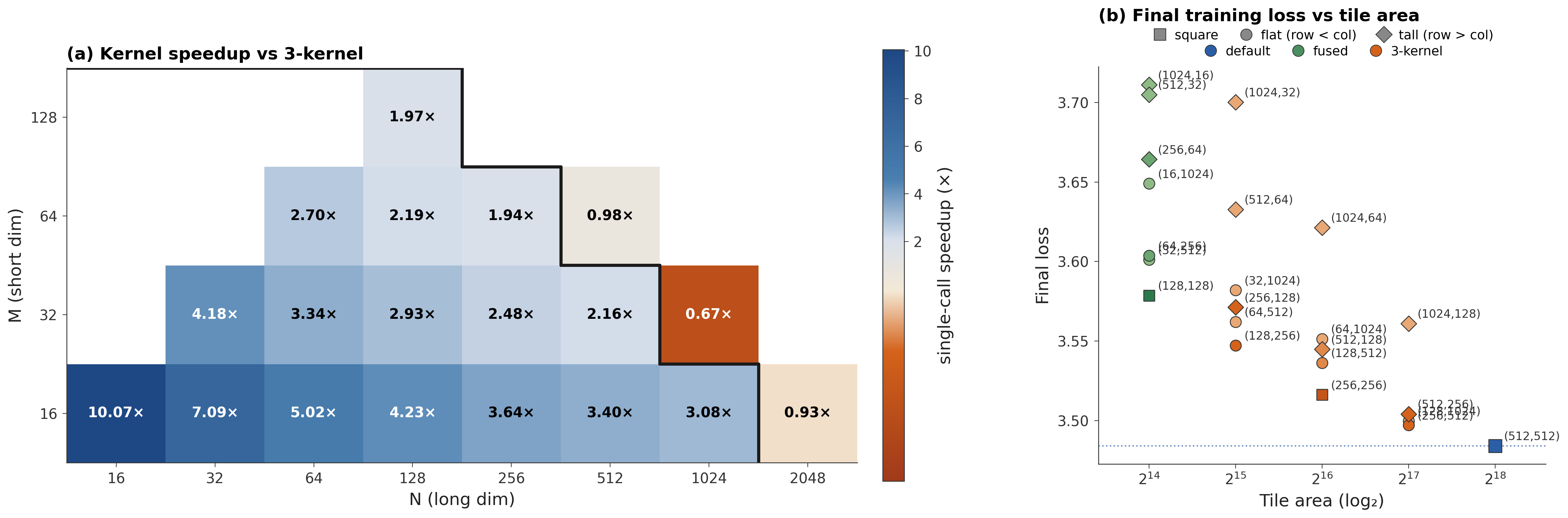}
    \caption{Speed and quality across tile shapes.
    \textbf{(a)} Single-call speedup of the fused NS5 kernel over the
    three-kernel path across shapes admissible under the A40 SRAM budget; the
    black outline marks the fused-dispatch region
    $(T_hT_w\le 16\,384)$.
    \textbf{(b)} Final training loss, averaged over the last $100$ of
    $300$ steps, for Llama-220M trained with HiMuon across tile shapes and
    kernel strategies.  The dotted line marks the $(512,512)$ reference.
    Shapes that are fastest in the isolated kernel are not necessarily best
    for training quality.}
    \label{fig:tile_shape_tradeoff}
\end{figure}

For the kernel-speed panel of \Cref{fig:tile_shape_tradeoff}, we sweep tile
shapes $(T_h,T_w)$ drawn from
$\{16,32,64,128,256,512,1024,2048\}$, keeping only shapes that fit the A40
SRAM budget.  The tile batch size is chosen so that the working set
$B_{\mathrm{tile}}T_hT_w$ is constant across shapes.  The reported speedup is
the fused NS5 kernel over the three-kernel path.  For the training-quality
panel, we train Llama-220M on FineWeb sample-10BT for $300$ steps at batch
size $8$, sequence length $1024$, and $\lr=0.02$, using a $20$-step linear
warmup followed by a $60$-step cosine decay to a $0.1$ learning-rate floor.
The experiment uses a single fixed seed and replays the same preloaded batch
sequence for every tile shape; the reported loss is the mean over the final
$100$ steps.

\Cref{fig:tile_shape_tradeoff} shows that kernel speed alone is not sufficient
to select a tile shape.  Smaller or more rectangular tiles can improve isolated
kernel speed, but may worsen training loss in this diagnostic.  In the tested
setting, square tiles give a more favorable quality--speed trade-off than
highly rectangular tiles.  We therefore use square $\tileSize=512$ tiles as the
default on A40.

\subsection{Cross-Layer Batching}
\label{sec:exp_xlayer}

We next evaluate the cross-layer batching strategy described in
\Cref{sec:xlayer}.  This experiment separates the gain from tiling itself from
the additional gain obtained by aggregating many independent tile computations
into larger GPU launches.

In Qwen3-1.7B, there are $28$ transformer blocks and $7$ Muon-eligible
matrices per block, giving $196$ independent per-layer Newton--Schulz calls in
an unbatched implementation.  With $\tileSize=512$, each individual matrix
contains only a modest number of tiles, so launching one Newton--Schulz call
per matrix can underutilize the GPU.  Cross-layer batching groups compatible
tiles across layers and processes them in a small number of larger batched
calls.

Tiles of the same shape and dtype are grouped across layers and concatenated
into a batched buffer.  When the buffer exceeds the memory budget, it is split
into the fewest chunks that fit, with intermediate workspace reused across
steps.  At $\tileSize=512$, the $1680$ tiles of Qwen3-0.6B fit in a single
call, whereas the $5376$ tiles of Qwen3-1.7B are processed in two chunks of
$2688$.  In this benchmark, every optimizer is run on the same synthetic
workload: $8$ stacked bias-free \texttt{Linear($H$,$H$)} layers in bf16, fed
fresh random Gaussian gradients at each step.  HiMuon uses $\tileSize=512$ and
$K=5$.  Per-step time is the median over $50$ steps after $20$ warmup steps,
and peak memory is the maximum allocated over a short run.  The cross-layer
variant additionally captures the batched Newton--Schulz workspace in a CUDA
graph, which pools and reuses that workspace across steps.

\begin{figure}[htbp]
    \centering
    \includegraphics[width=0.85\textwidth]{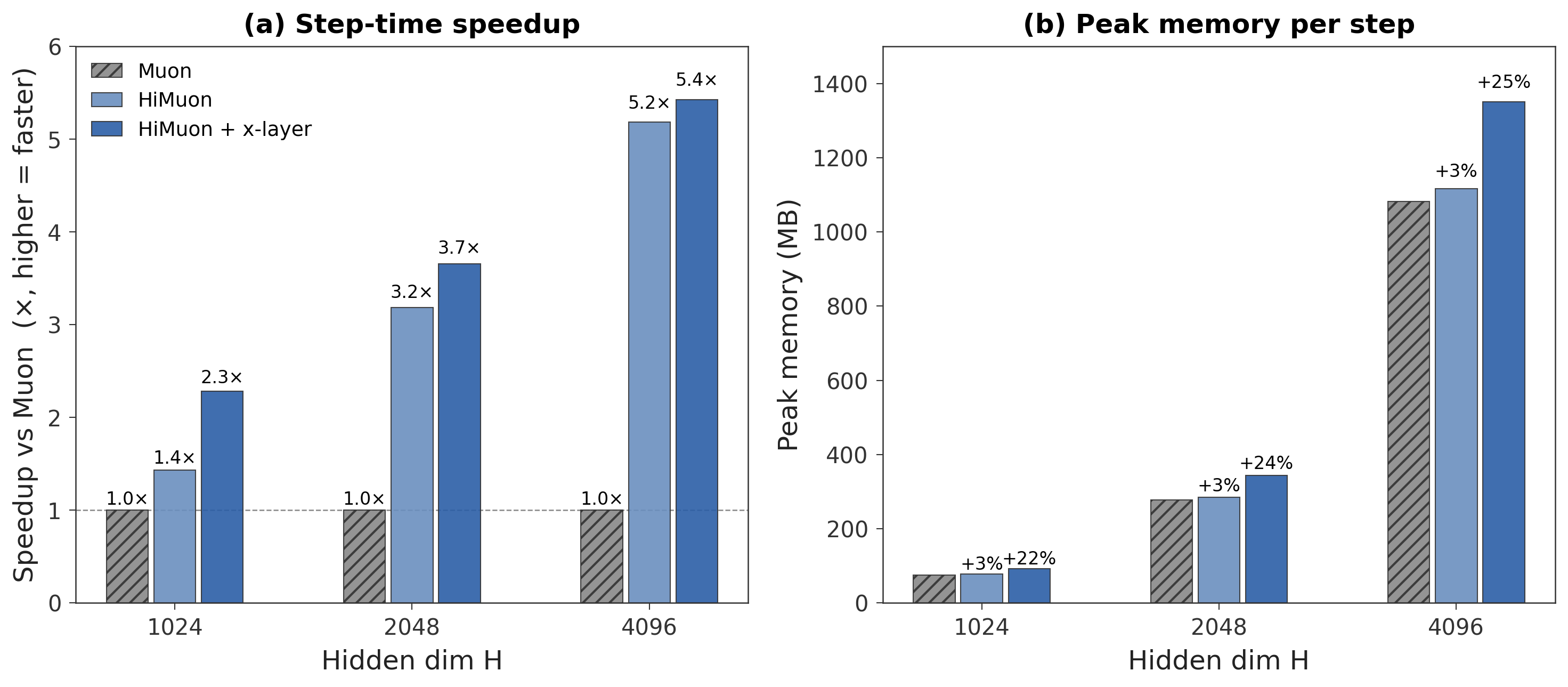}
    \caption{Optimizer-step performance on the NVIDIA A40 GPU.  The workload
    consists of $8$ stacked \texttt{Linear($H$,$H$)} layers in bfloat16, with
    $H\in\{1024,2048,4096\}$.  \textbf{(a)} Per-step speedup relative to
    full-matrix Muon.  Tiled HiMuon reaches
    $1.4\times$/$3.2\times$/$5.2\times$ at $H=1024/2048/4096$; cross-layer
    batching raises the total speedup to
    $2.3\times$/$3.7\times$/$5.4\times$.  \textbf{(b)} Peak GPU memory.
    Tiled HiMuon adds about $+3\%$ over Muon at all $H$; cross-layer batching
    adds about $+22\%$ to $+25\%$ due to the batched buffer, with the overhead
    bounded by chunking and workspace reuse.}
    \label{fig:optimizer_step_bench}
\end{figure}

\Cref{fig:optimizer_step_bench} shows that the tiled map provides the dominant
optimizer-step speedup.  Cross-layer batching gives an additional wall-clock
reduction by replacing many small Newton--Schulz calls with fewer larger
batched calls.  The incremental benefit is largest at $H=1024$, where the
plain tiled calls are small; it shrinks at $H=4096$, where each call is already
large enough to better utilize the GPU.  The batching buffer increases peak
memory, but the overhead is controlled by chunking and workspace reuse and does
not grow with the number of training steps.

\subsection{Capacity Robustness to Low Precision}
\label{sec:precision_recovery}

Low-precision arithmetic in the Newton--Schulz inner GEMM is attractive on
modern accelerators, but it perturbs the matrix-function update through
rounding error~\cite{HighamMary2019}.  We assess
whether this perturbation changes the optimizer's behavior on a controlled
downstream metric.  Using the one-step associative-memory benchmark of
\Cref{sec:capacity}, we hold the problem instance and gradient fixed within
each seed and vary only the GEMM precision used inside the Newton--Schulz
iteration.  Thus differences in recovered items reflect finite-precision
arithmetic in the map rather than a different problem instance or tuned
learning rate.

The setup uses $d=1024$, $N=4096$ keys, Zipf exponent $\alpha=1.5$,
batch-to-dimension ratio $B/d=10$, $K=5$ Newton--Schulz iterations, and five
random seeds, executed on a single NVIDIA L40S.  For each seed, the
key--value pairs
$(v_i,u_i)\sim\mathcal{N}(0,d^{-1}\mI_d)$ and the resulting minibatch gradient
$\mG_0$ are constructed in fp64 and held fixed across precision variants.  The
four precision settings are an fp64 reference, fp32, bf16, and fp8e4m3 with
per-tensor scaling.  The learning rate $\eta=1.0$ is chosen once from the
fp64 reference run and reused for every precision.  We compare full-matrix
Muon and HiMuon with $\tileSize\in\{128,256,512\}$.

\begin{table}[htbp]
  \centering
  \footnotesize
  \caption{One-step recovered items on the associative-memory benchmark of
  \Cref{sec:capacity} at $d=1024$, $N=4096$, $\alpha=1.5$, $B/d=10$,
  $K=5$, five random seeds, NVIDIA L40S.  Recovered-item entries are
  mean~$\pm$~seed standard deviation.  The last two columns give the operator
  gap
  $\delta_F =
  \|\mW_1^{(\cdot)}-\mW_1^{\mathrm{fp64}}\|_F/
  \|\mW_1^{\mathrm{fp64}}\|_F$,
  averaged over seeds, between the produced one-step update and the fp64
  reference of the same update rule.  Problem embeddings, minibatch indices,
  and gradient $\mG_0$ are shared across precision columns within each seed;
  only the GEMM precision inside the Newton--Schulz iteration differs.}
  \label{tab:precision_recovery}
  \begin{tabular}{lcccc cc}
    \toprule
    & \multicolumn{4}{c}{Recovered items}
    & \multicolumn{2}{c}{Op.\ gap $\delta_F$ to fp64} \\
    \cmidrule(lr){2-5} \cmidrule(lr){6-7}
    Method        & fp64                & fp32                & bf16                & fp8e4m3             & bf16  & fp8e4m3 \\
    \midrule
    Muon          & $517.8\!\pm\!11.3$  & $517.8\!\pm\!11.3$  & $517.4\!\pm\!11.7$  & $517.2\!\pm\!11.4$  & 0.062 & 0.019 \\
    HiMuon-$512$  & $508.8\!\pm\!12.1$  & $508.8\!\pm\!12.1$  & $508.2\!\pm\!11.8$  & $509.0\!\pm\!11.9$  & 0.061 & 0.027 \\
    HiMuon-$256$  & $486.4\!\pm\!12.6$  & $486.4\!\pm\!12.6$  & $486.4\!\pm\!11.2$  & $485.8\!\pm\!12.6$  & 0.058 & 0.039 \\
    HiMuon-$128$  & $408.0\!\pm\!15.3$  & $408.0\!\pm\!15.3$  & $406.6\!\pm\!15.6$  & $406.8\!\pm\!14.6$  & 0.052 & 0.059 \\
    \bottomrule
  \end{tabular}
\end{table}

\Cref{tab:precision_recovery} shows that, for every update rule, the recovered
item counts agree across precision settings to well within the seed-to-seed
standard deviation.  fp32 reproduces the fp64 means exactly under this setup.
bf16 and fp8e4m3 produce operator-level gaps of roughly $2\%$--$6\%$ relative
to the per-method fp64 reference, yet change the mean recovered count by at
most about one item out of roughly $400$--$520$ recovered items.  This is much
smaller than the seed-to-seed standard deviation of roughly $11$--$16$ items.
Thus, in this diagnostic, using bf16 or fp8e4m3 for the Newton--Schulz inner
GEMMs does not measurably change the downstream recovery metric.

\subsection{Alignment and Descent Potential under Low Precision}
\label{sec:precision_descent}

The recovery test in \Cref{sec:precision_recovery} examines one downstream
metric.  We also evaluate two local quantities associated with per-step
progress of Muon-type updates: the alignment $\gamma$ of the update with the
true gradient and the descent potential $\phi$ along the update direction, as
studied by Shumaylov et al.~\cite{ShumaylovMuonNotSpecial2026}.  This
experiment asks whether tile size and low-precision Newton--Schulz arithmetic
change these quantities on a controlled quadratic proxy.

We write
$\langle\mathbf{A},\mathbf{B}\rangle_F=\tr(\mathbf{A}^\top\mathbf{B})$ for the
Frobenius inner product.  Given a full-batch gradient $\mG$, a minibatch
gradient $\widetilde{\mG}$, an update direction $\mathbf{D}$, and Hessian
action $\mathcal{H}[\mathbf{D}]$, define
\[
  \gamma
  =
  \frac{\langle\mG,\mathbf{D}\rangle_F}
       {\langle\widetilde{\mG},\mathbf{D}\rangle_F},
  \qquad
  \phi
  =
  \frac{\langle\widetilde{\mG},\mathbf{D}\rangle_F^2}
       {\langle\mathbf{D},\mathcal{H}[\mathbf{D}]\rangle_F}.
\]
Larger values correspond to faster local progress under an appropriately tuned
step size.  We evaluate these quantities on a teacher--student quadratic
least-squares proxy after Davis and Drusvyatskiy~\cite{P1243}:
\[
    f(\mW)
    =
    \frac12
    \mathbb{E}_{\boldsymbol{x}\sim\mathcal{N}(0,\boldsymbol{\Sigma})}
    \left\|(\mW-\mW^\star)\boldsymbol{x}\right\|^2 .
\]
The anisotropic covariance $\boldsymbol{\Sigma}$ and teacher $\mW^\star$ yield
the full gradient, minibatch gradient, and Hessian action in closed form.

\noindent\textbf{Protocol.}
We use $d=1024$, spectrum $\lambda_i\propto i^{-\alpha}$ with $\alpha=1$ and
unit mean, a random orthogonal eigenbasis, minibatch size $B=256$, and one
seed.  The teacher $\mW^\star$ and initial iterate are random unit-Frobenius
matrices.  A single fp32 full-matrix Muon trajectory of $40$ steps
($\lr=0.03$, unit-Frobenius-normalized update) generates the iterates.  At
each step, every tile-size and precision variant sees the same minibatch
gradient $\widetilde{\mG}$.  The evaluated update is
$\mathbf{D}=\Phi_{K,\tileSize}(\widetilde{\mG})$ for
$\tileSize\in\{128,256,512\}$, or the full map
$\Phi_K(\widetilde{\mG})$ on the full $1024\times1024$ matrix, with $K=5$.
The map is evaluated in fp32, bf16, and fp8e4m3 with per-tile dynamic scaling
on an NVIDIA L40S.

\begin{figure}[htbp]
    \centering
    \includegraphics[width=0.85\textwidth]{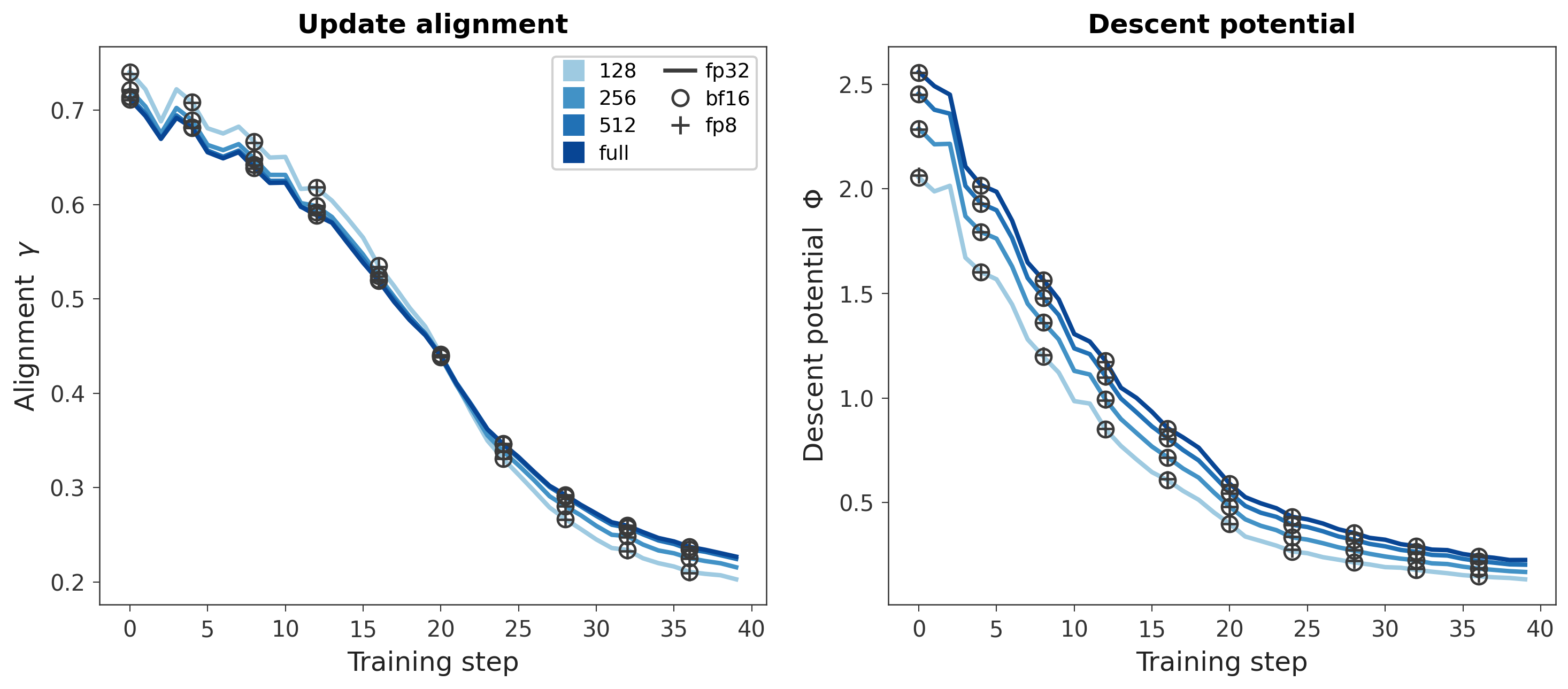}
    \caption{Update alignment $\gamma$ \textbf{(left)} and descent potential
    $\phi$ \textbf{(right)} along a $40$-step fp32 full-matrix Muon trajectory
    on a teacher--student quadratic proxy after Davis and
    Drusvyatskiy~\cite{P1243}.  The setup uses $d=1024$, $\alpha=1$,
    $B=256$, $K=5$, one seed, and an NVIDIA L40S.  Line color encodes the
    Newton--Schulz operator, from $\tileSize=128$ to $\tileSize=512$ and the
    full $1024\times1024$ map.  Precision is overlaid as bf16 circles and fp8
    crosses on the fp32 line, sampled every fourth step.}
    \label{fig:precision_descent}
\end{figure}

\Cref{fig:precision_descent} shows that both $\gamma$ and $\phi$ decline along
the trajectory as the iterate approaches the teacher.  The precision variants
are visually indistinguishable at the plotted resolution: bf16 and fp8 markers
lie on the corresponding fp32 curves for each tile size.  Tile size has a more
visible effect.  The alignment curves remain close across tile sizes, while
the descent-potential curves follow the ordering
$\text{full}>512>256>128$, with the largest gap early in the trajectory.  This
matches the qualitative ordering seen in the capacity diagnostic of
\Cref{sec:capacity_diagnostic}.  These results suggest that, in this controlled
proxy, tile size changes the local update geometry more than the tested
Newton--Schulz GEMM precisions.

\section{Conclusion}
\label{sec:conclusion}

In this paper, we introduced HiMuon, a tiled Newton--Schulz update for
Muon-type optimization.  HiMuon replaces one full-matrix Newton--Schulz map by
independent applications of the same finite map on fixed-size intra-matrix
tiles.  This defines a local matrix-function update whose tile size controls
the trade-off between retained spectral coupling and computational cost.  For
fixed $\tileSize$, the leading Newton--Schulz work scales as
$\mathcal{O}(HW\tileSize K)$ rather than
$\mathcal{O}(HW\min\{H,W\}K)$.

The experiments characterize this trade-off at several levels.  Capacity
diagnostics show that larger tiles more closely preserve the behavior of
full-matrix Muon, while smaller tiles reduce cost.  GPU experiments show that
the tiled formulation can be implemented efficiently through
tile-size-dependent kernels, cross-layer batching, and memory-bounded
chunking.  End-to-end Qwen3 training experiments indicate that HiMuon reduces
optimizer-step and wall-clock time while maintaining training and validation
losses close to full-matrix Muon in the tested single-seed runs.

Future work includes a theory of the tiled map, adaptive schedules that choose
$\tileSize$ from training state or layer shape, and broader evaluation across
model families, hardware platforms, and distributed training layouts.

\section*{Acknowledgments}
The authors acknowledge the Minnesota Supercomputing Institute (MSI) at the University of Minnesota for providing the computational resources on which the majority of the experiments in this paper were run.
The authors thank the developers of the Muon optimizer for making their implementation publicly available.

\bibliographystyle{siamplain}
\bibliography{papers}

\end{document}